\documentclass[12pt]{article}
\usepackage{amssymb,amsmath, txfonts}
\usepackage{appendix}
\begin{document}

\def\abstractname{\bf Abstract}
\def\dfrac{\displaystyle\frac}
\let\oldsection\section
\renewcommand\section{\setcounter{equation}{0}\oldsection}
\renewcommand\thesection{\arabic{section}}
\renewcommand\theequation{\thesection.\arabic{equation}}
\newtheorem{theorem}{\indent Theorem}[section]
\newtheorem{lemma}{\indent Lemma}[section]
\newtheorem{proposition}{\indent Proposition}[section]
\newtheorem{definition}{\indent Definition}[section]
\newtheorem{remark}{\indent Remark}[section]
\newtheorem{corollary}{\indent Corollary}[section]
\def\pd#1#2{\displaystyle\frac{\partial#1}{\partial#2}}
\def\d#1{\displaystyle\frac{d#1}{dt}}

\title{\LARGE\bf Long time dynamics of the Cauchy problem for the predator-prey model with cross-diffusion
\\
\author{Chunhua Jin$^a$, Yifu Wang$^b$
\thanks{
Corresponding author. Email: {\tt wangyifu@bit.edu.cn}}
\\
\small \it{$^a$School of Mathematical Sciences, South China
Normal University, }
\\
\small \it{Guangzhou, 510631, P. R. of China}
\\
\small \it{$^b$
School of Mathematics and Statistics, Beijing Institute of Technology, }
\\
\small \it{Beijing 100081,  P. R. of China}
}}

\date{}

\maketitle

\begin{abstract}
This paper is concerned with a predator-prey model in $N$-dimensional spaces ($N=1, 2, 3$), given by
\begin{align*}\left\{\begin{aligned}
&\frac{\partial u}{\partial t}=\Delta u-\chi\nabla\cdot(u\nabla v),\\
&\frac{\partial v}{\partial t}=\Delta v+\xi\nabla\cdot(v\nabla u),
\end{aligned}\right. \end{align*}
which describes  random movement of both predator and prey species, as well as the spatial dynamics involving predators pursuing prey and prey attempting to evade predators. It is shown that any global strong solutions of the corresponding Cauchy problem converge to zero in the sense of $L^p$-norm for any $1<p\le \infty$, and also converge to the heat kernel with respect to $L^p$-norm  for any $1\le p\le \infty$. In particular,  the decay rate thereof is optimal in the sense that
it is  consistent with that of the heat equation in $\mathbb R^N$ ($N=2, 3$).
Undoubtedly, the global existence of solutions  appears to be among  the most challenging topic in the analysis of this model.
 Indeed even in the one-dimensional setting,  only global weak solutions in a bounded domain have been successfully constructed by far. Nevertheless, to provide a comprehensive understanding of the main results, we append the conclusion on the global existence and asymptotic behavior of strong solutions,  although  certain smallness conditions on the initial data are required.
\end{abstract}



\section{Introduction}
The interaction between different species generates a complex ecosystem, and the predator-prey relation
  is one of the fundamental binary interactions, in which predators pursue preys. Since Lotka-Volterra-type prey-predator system was introduced,
   various mathematical models have
been proposed to explore a diverse range of complex ecological phenomena, including Allee effects \cite{Stephens}, group defense \cite{Venturino} and so on, and especially  predator-prey model invoking some mechanism of spatial interaction
is derived to describe complex space-time patterns which reflect nonlocal spatial interactions between the predators and preys
 (see \cite{Okubo, NWM, Murray} for example). Therefore the investigation of the diffusive predator-prey systems becomes an  important step
through which one can get a rather deeper understanding of the effects of the spatial  dispersal of the predators and preys on
 the entire ecosystem.
In the spatial  predator-prey relationship,
 besides the kinetic dynamics and the random movement of the predators and preys,
  the predators and preys may also partially orient
their movement towards or away from higher concentrations of the other species, that is predators move toward the gradient direction of prey distribution (called ``prey-taxis"), and/or preys move opposite to the gradient of predator distribution (called ``predator-taxis"). To capture this pursuit-escape phenomenon,
Tsyganov et.al \cite{TBH} proposed the cross-diffusion system of the form
\begin{align}\label{TS}
\left\{\begin{aligned}
&\frac{\partial u}{\partial t}=\Delta u-\chi\nabla\cdot(u\nabla v)+T_1(u, v),\\
&\frac{\partial v}{\partial t}=\Delta v+\xi\nabla\cdot(v\nabla u)+T_2(u, v),
\end{aligned}\right.\end{align}
where the unknown functions $u(t, x)$ and $ v(t, x) $ correspond to the predator and prey population densities respectively at position $x$ and time $t>0$;
$ -\chi\nabla\cdot(u\nabla v)$ and $+\xi\nabla\cdot(v\nabla u)$
 stand for the predators are attracted by the preys with coefficient $\chi$ and the preys escape from predators with coefficient $\xi$, respectively.  The
functions $T_1(u, v)$ and $T_2(u, v) $ represent the inter-specific interaction of $u$ and $v$.

As pointed in \cite{TW1}, \eqref{TS} can be regarded as a distant relative of the classic Shigesada-Kawasaki-Teramoto  model (\cite{SKT})
\begin{align*}
\left\{\begin{aligned}
&\frac{\partial u}{\partial t}=\Delta[(a_1+a_2u+a_3v)u]+T_1(u, v),
\\
&\frac{\partial v}{\partial t}=\Delta[(b_1+b_2u+b_3v)v]+T_2(u, v),
\end{aligned}\right.
\end{align*}
abbreviated as the SKT model. The SKT model and its simplified version have been widely studied in the past few decades, see for example \cite{DP, KJ, LN, LWY, LWI} or the references therein. However up to now, it seems that the global well-posedness theory of solutions thereof has not been fully established yet.
Unlike the SKT model, the cross-diffusion in \eqref{TS} originates from a specific  mechanism  which characters the attractive and repulsive taxis between the two populations simultaneously. As a remark feature, \eqref{TS} possesses
a quite colorful wave-like solution,  for example, 
quasi-soliton pattern is identified in \cite{TBH}.
In recent years, the global well-posedness theory of
 \eqref{TS} has received a considerable attention from researchers.
However,  the achievement on  system \eqref{TS} seems
at a rather rudimentary stage with rather limited results available for the references \cite{FM, TW1, TW2, XL1, XL2}.  For instance, for the case  $T_1=T_2=0$ or $T_1$ and $T_2$ being  Lotka-Volterra response terms, Fuest investigated the existence of global classical solutions of homogeneous Neumann boundary value problem of
 \eqref{TS} on a smooth bounded domain   when the initial values are sufficiently close to the constant equilibrium points, as well as the asymptotic behavior thereof in \cite{FM}.
 Nevertheless, due to the double tactic mechanism in  \eqref{TS},   there is significant challenge in the rigorous analysis of \eqref{TS} with large initial data. Indeed, to the best of our knowledge, the global solvability for the initial boundary value problem  for \eqref{TS} with large initial data is limited in one-dimensional spatial situation yet, see \cite{TW1, TW2, XL1} for example.

In contrast to the sparse results on model \eqref{TS} by far, 
extensive studies have been carried out on  the simplified variants thereof.
As an important particular case of \eqref{TS} in which the predator-taxis is neglected,  the  diffusive predator-prey system only with prey-taxis reads as
 \begin{align}\label{1.2}
\left\{\begin{aligned}
&\frac{\partial u}{\partial t}=\Delta u-\chi\nabla\cdot(u\nabla v)+T_1(u, v),
\\
&\frac{\partial v}{\partial t}=\Delta v+T_2(u, v).
\end{aligned}\right.\end{align}
Here the most type of  ecological inter-specific interaction can be encapsulated by the choice of $T_1(u,v)=\gamma  uF(v)-uh(u)$,  $T_2(u,v)=- uF(v) +k(v) $, $F(v)$ is the predation rate, the parameter $\gamma>0$ denotes the conversion rate from prey to predator, $uh(u)$ and $k(v)$ are referred to the intra-specific interaction of  the predators and preys, respectively.
Up to date there are a few results available to  \eqref{1.2} in the bounded domain 
(see \cite{CA,He,Jinwang,Wang,WinklerJDE,Wu} for example and references cited therein). For instance, the global existence and stability of corresponding solutions to  \eqref{1.2} with no-flux boundary conditions was established in a two-dimensional bounded domain under some specific condition on $F,h$ and $k$ (\cite{Jinwang}), while the existence of global bounded solutions was obtained in higher dimensions for small $\chi>0$ and in three-dimensional setting for suitably small initial value $\|v_0\|_{L^\infty(\Omega)}$, respectively (see \cite{LiDan,Wu}). In particular, Winkler (\cite{WinklerJDE}) showed the global existence of weak solutions over a bounded domain with space dimension $N\leq 5$ and their eventual smoothness under some weak conditions which includes $k(v)=0$ and $h(u)=0$.

On the other hand, prey-taxis model \eqref{1.2} can also be regarded as the extension of  the classical Keller-Segel model (\cite{KS})
\begin{align}\label{1.3}
\left\{
\begin{aligned}
&\frac{\partial u}{\partial t}=\Delta u-\chi\nabla\cdot(u\nabla v),&\quad x\in \Omega,t>0,
\\
&\frac{\partial v}{\partial t}=\Delta v-v+u, &\quad x\in \Omega,t>0,
\end{aligned}\right.
\end{align}
where $\Omega$ is either a bounded domain  in  $\mathbb{R}^N$ with smooth boundary $\partial\Omega$ or the whole space $\mathbb{R}^N$, which
describes a chemotactic feature, the aggregation of some organisms (cellular slime molds) denoted by $u$
sensitive to the gradient of a chemical substance denoted by $v$. The possibly most striking feature of this model is
the occurrence of a critical mass phenomenon in two dimensions, namely there exists  a threshold value $M_c>0$ with respect to
initial datum $\|u_0\|_{L^1}$  such that  all solutions exist globally in time for
$\|u_0\|_{L^1}\leq M_c$, while  finite time blow-up occurs if  $\|u_0\|_{L^1}>M_c$. Here
$M_c$ is either $\frac{8\pi}{\chi}$ or $\frac{4\pi}{\chi}$ according to $\Omega$ is  the whole space or the bounded domain
\cite{Horstmann,NSY,Mizoguchi}. Notice that the approach to dealing with the dynamical behavior of  \eqref{1.3} in the whole space $\mathbb{R}^N$ $(N\geq 3)$ is slightly different from that on the bounded domain. For instance, it is well known
that  the functional
$$
\mathcal{F}(u,v):=\int u\log u  - \chi \int u v +\frac \chi 2 \int  v^2+  \frac \chi 2 \int  |\nabla v|^2
$$
satisfies
\begin{align} \label{1.4}
\frac {d}{dt}\mathcal{F}(u,v)=
- \mathcal{D}(u,v):= -\int | \triangle v - v + u|^2 - \int \left|\frac {\nabla u}{\sqrt{u}}- \chi \sqrt{u} \nabla v\right|^2
\end{align}
along reasonably regular solution curves of \eqref{1.3}, which then allows us for identifying the circumstances under which solutions exist globally and remain bounded. Going beyond that, \eqref{1.4} can also serve as
a fundamental ingredient in the detection of blow-up phenomena. Indeed, if there are some $C>0$ and $\vartheta \in (0, 1)$ such that
 \begin{align}\label{1.5}
   \mathcal{F} (u,v)\geq -C \cdot ( \mathcal{D}^\vartheta(u,v) + 1)
\end{align}
 throughout some radially symmetric functions $(u,v)$, \eqref{1.4} is then capable of identifying large
sets of initial data which lead to finite-time blow-up in Neumann problems for \eqref{1.3} (\cite{WinklerJMPA}).
However, as point out in \cite{WinklerNonlinearity}, term $\int_{\mathbb{R}^2} u\log u$ in $\mathcal{F}(u,v)$
may be no longer bounded from below,  hence the energy argument above is not accessible  for the Cauchy problem.

The investigation of long-time dynamic properties in biological models is crucial for comprehending the fundamental processes of population growth, death, and migration within biological populations. Additionally, it sheds light on how these processes impact the stability and structure of ecosystems, this research also serves as a valuable resource for ecological conservation and management practices.
Regarding the global solution of this kind of model on a bounded domain, the homogeneous Neumann boundary value problem has been thoroughly investigated as we stated above. It is well-known that the asymptotic behavior  under zero Neumann
boundary condition is nearly the same as the corresponding ODE systems, that is the solution  will tend to some constant steady
state eventually. In such cases,  numerous studies have demonstrated, see for example \cite{Jinwang, WM1} or the references therein. In addition, some necessary conditions for pattern formation have also been explored \cite{KOS, LHL}. However, the long-time  asymptotic behavior of the Cauchy problem for solutions in the whole space differs significantly from that observed on bounded domains.
As an illustrative example, considering the 2D Patlak-Keller-Segel system
\begin{align*}
\left\{
\begin{aligned}
&\frac{\partial u}{\partial t}=\Delta u-\chi\nabla\cdot(u\nabla v),
\\
&\Delta v+u=0,
\end{aligned}\right.
\end{align*}
an important observation in this model is its scaling invariance property, meaning that if $(u(x,t), v(x,t))$ is a solution, so is
$(u_\lambda(x,t), v_{\lambda}(x,t))$, where
$$
u_\lambda(x,t)=\lambda^{-2}u(\lambda^{-1}x, \lambda^{-2}t), \quad v_{\lambda}(x,t)=v(\lambda^{-1}x, \lambda^{-2}t).
$$
This is a crucial property that significantly influences its dynamic behavior.
Biler et al. \cite{BKL, BM1, BP, MM, MN} have delved into the existence and characteristics of self-similar solutions in the form of $t^{-1}G_\alpha\left(\frac{|x|}{\sqrt t}\right)$ when the initial mass $\alpha$ is below the critical threshold in two dimensional space. It has been confirmed that such self-similar solutions also act as global attractors \cite{BDP, BM1}, that is
$$
\lim_{t\to\infty}\left\|u(x,t)-t^{-1}G_\alpha\left(\frac{|x|}{\sqrt t}\right)\right\|_{L^p}=0
$$
for $p\in [1,\infty]$ if $\|u\|_{L^1}=\alpha$.
Furthermore, the optimal rate of convergence towards this self-similar profile has been presented \cite{BM, FM}.
Similarly, for the following parabolic-parabolic Keller-Segel system,
\begin{align*}
\left\{
\begin{aligned}
&\frac{\partial u}{\partial t}=\Delta u-\chi\nabla\cdot(u\nabla v)
\\
&\frac{\partial v}{\partial t}=\Delta v+u,
\end{aligned}\right.
\end{align*}
it is shown that the solution will finally converge to the self-similar solution of this system for small mass in $\mathbb R^2$ \cite{KM, NY}. Similar results have also been observed in the two-dimensional Navier-Stokes equations \cite{GW}, and  some reaction-diffusion equations \cite{EKM, GV, KP}.
While for the following Keller-Segel system
\begin{align*}
\left\{
\begin{aligned}
&\frac{\partial u}{\partial t}=\Delta u-\chi\nabla\cdot(u\nabla v),
\\
&\frac{\partial v}{\partial t}=\Delta v-v+u,
\end{aligned}\right.
\end{align*}
it is evident that this equation does not exhibit scaling invariance, hence no self-similar solutions exist. Regarding the long-time asymptotic behavior of this model,   Nagai \cite{NA} in 2001  investigated the asymptotic profile of bounded solutions
in $\mathbb R^2$, and demonstrated that the solution converges to the Gauss kernel in terms of the $L^\infty$-norm. Later, in 2003, Nagai, Syukuinn et.al \cite{NSU} extended this result to $\mathbb R^N$ with $N\ge 2$, demonstrating convergence in $L^p$-norm for $1<p\le \infty$.
Recently, Carrillo et al. \cite{CCY} studied a general aggregation-diffusion equation of the form
$$
\frac{\partial u}{\partial t}=\Delta u+\nabla\cdot(u\nabla W*u),
$$
and proved that if the interaction potential is bounded and its first and second derivatives decay fast enough at infinity, then the linear diffusion eventually dominates its attractive or repulsive effects, and thereby solutions finally behave like the fundamental solution of the heat equation.

Motivated by the above works,  this paper considers the long-time dynamical behavior of strong solutions to the Cauchy problem of \eqref{TS} in the higher dimensions, which reads as
\begin{align}\label{1.6}
\left\{
\begin{aligned}
&\frac{\partial u}{\partial t}=\Delta u-\chi\nabla\cdot(u\nabla v) && \text{in}\  Q,
\\
&\frac{\partial v}{\partial t}=\Delta v+\xi\nabla\cdot(v\nabla u) && \text{in}\  Q,
\\
& u(x,0)=u_0(x),  v(x,0)=v_0(x) && \text{in}\ \mathbb R^N,
\end{aligned}\right.
\end{align}
where $Q=\mathbb R^N\times \mathbb R^+$ with $N=1, 2, 3$, and the initial value $u_0, v_0\in L^1(\mathbb R^N)\cap L^\infty(\mathbb R^N)$ are nonnegative. Our focus is that in the whole space, how far the attractive and repulsive interaction of  the species $u$ and $v$
 affect their long dynamical behavior in the situation that their populations dynamics  is neglected.

 It is observed that the system \eqref{1.6} adheres to the law of conservation of mass, that is
\begin{align}\label{1.7}
\int_{\mathbb R^N}u(x,t)dx=\int_{\mathbb R^N}u_0(x)dx,\quad \int_{\mathbb R^N}v(x,t)dx=\int_{\mathbb R^N}v_0(x)dx
\end{align}
for all $t>0$, whereas our findings reveals  that the bounded strong solutions of this system converge to zero in the sense of $L^p$-norm with any $1<p\le \infty$. Apart from that,  due to the fact that the $L^p$-norm of the heat kernel $G(\cdot, t)$ decreases to $0$ as $t\to\infty$ for any $p>1$, that is
$$
\|G(\cdot, t)\|_{L^p}=(4\pi t)^{-\frac N2}\left(\int_{\mathbb R^N}\exp\left(-\frac{p|x|^2}{4t}\right)dx\right)^{\frac1p}
= \pi^{-\frac N2}(4t)^{-\frac N2(1-\frac1p)}\left(\int_{\mathbb R^N}\exp(-py^2)dy\right)^{\frac1p}\to 0,
$$
 the strong solutions hence converges to the heat kernel in $L^p$-space for any $p>1$.
Therefore a natural question is whether this solution also converges to the heat kernel in the sense of $L^1$-norm.
To this end, the more rigorous analysis is required. Indeed, we not only prove that
for any $p$ with $1\le p\le \infty$, the solutions converge to the heat kernel in $L^p$-norm sense, but also provide  accurate convergence rate estimates and  the decay rate,  which appears to be optimal in the sense that the rate achieved herein is  consistent with that of the heat equation in $\mathbb R^N$ ($N=2,3$) except for $N=1$.
 In addition,  as the complement to the above results,  it is shown  that the global existence of strong solutions to \eqref{1.6} in the space of $L^\infty(\mathbb R^+; H^2(\mathbb R^N))$ for some small initial data. Here, we define a strong solution as a pair $(u, v)$ belonging to $L^\infty(\mathbb R^+; H^2(\mathbb R^N))$, which not only satisfies the system of equations \eqref{1.6} in the sense of distributions,  but also meets the  equations almost everywhere.

In Section 2, we employ a scaling transformation  and  $L^p-L^q$ estimates for $e^{t\Delta}$ to derive the decay estimate of any strong solution $(u, v)\in L^\infty(\mathbb R^+; H^2(\mathbb R^N)\cap L^1(\mathbb R^N))$ of \eqref{1.6}, and thereby obtain the following convergence results by means of a bootstrap argument.
\begin{theorem}
\label{thm-1}
Let $N=2, 3$ and  $(u, v)\in L^\infty(\mathbb R^+; H^2(\mathbb R^N)\cap L^1(\mathbb R^N))$ be the global strong solution of \eqref{1.6}.
Then there exists $C>0$ such that
\begin{align*} \begin{aligned}
&\|u(t)\|_{L^p}\le  C(t+1)^{-\frac{N}2\left(1-\frac1p\right)}\|u_0\|_{L^1},
\\
&\|v(t)\|_{L^p}\le  C(t+1)^{-\frac{N}2\left(1-\frac1p\right)}\|v_0\|_{L^1}\end{aligned}
\end{align*}
where $p\in (1,\infty]$.

Moreover, for any $p\in [1,\infty]$
\begin{align*} \begin{aligned}
&\lim_{t\to\infty}t^{\frac N{2}\left(1-\frac1p\right)}\left\|u(\cdot,t)-\int_{\mathbb R^N}u_0(y)dy G(\cdot, t)\right\|_{L^p}=0,
\\
&\lim_{t\to\infty}t^{\frac N{2}\left(1-\frac1p\right)}\left\|v(\cdot,t)-\int_{\mathbb R^N}v_0(y)dy G(\cdot, t)\right\|_{L^p}=0,
\end{aligned}
\end{align*}
where $G(x,t)$ is the heat kernel, i.e. $
G(x,t)=(4\pi t)^{-\frac N2}\exp\left(-\frac{|x|^2}{4t}\right).
$
\end{theorem}

In the above theorem, we obtained an estimate of the optimal decay rate for the solutions in two and three dimensions. However, for the one-dimensional spatial case, the bootstrap technique  cannot be employed to achieve an optimal decay rate, and only the following conclusion has been reached.

\begin{theorem}
\label{thm-2}
Let $N=1$. Suppose that $(u, v)
\in L^\infty(\mathbb R^+; H^2({\mathbb R} )\cap L^1(\mathbb R ))$ is  a global strong solution of \eqref{1.6}.
Then there exists $C>0$, such that
\begin{align*}\begin{aligned}
&\|u(t)\|_{L^p}\le  C(t+1)^{-\frac{3}8\left(1-\frac1p\right)}\|u_0)\|_{L^1},
\\
&\|v(t)\|_{L^p}\le  C(t+1)^{-\frac{3}8\left(1-\frac1p\right)}\|v_0\|_{L^1},
\end{aligned}\qquad
\end{align*}
where $p\in (1,\infty]$. Additionally, for any $0<m<\frac{9}{32}$
\begin{align*}
\lim_{t\to\infty}t^{m\left(1-\frac1p\right)}\left\|u(\cdot,t)-\int_{\mathbb R^N}u_0(y)dy G(\cdot, t)\right\|_{L^p}=0,
\end{align*}
where $p\in (1,\infty]$.
\end{theorem}

As a necessary complement to the above result, we make sure that the strong solution of \eqref{1.6}
in the space of $L^\infty(\mathbb R^+; H^2(\mathbb R^N))$ actually exists in Section 3.

\begin{theorem}
\label{thm-3}
Assume that $N=1, 2, 3$,  the nonnegative initial value $(u_0, v_0)$ satisfies $u_0, v_0\in L^1(\mathbb R^N)\cap H^2(\mathbb R^N)$ and $|x|^ru_0(x)\in L^2(\mathbb R^N)$ with some constant $r\in (0, 1)$. Denote
$$
\phi(t)=\|u(\cdot, t)\|_{H^2}^{2}+\|v(\cdot, t)\|_{H^2}^{2},
$$
$$
h(t)=\|\nabla u(\cdot,t)\|_{H^2}^2+\|\nabla v(\cdot,t)\|_{H^2}^2.
$$
Then if $\phi(0)\le \frac1{2C^*(\chi^2+\xi^2)}$ for some constant $C^*$ depending only on $N$, the problem \eqref{1.6} admits a unique  global strong solution $(u, v)$ with $u, v\in \mathcal X$ such that
$$ \quad \frac12 \int_0^\infty h(t)dt\le \phi(0),\quad
\phi(t)\le \left(\phi^{-\frac 2N}(0)+\frac{2}{N}A^*t\right)^{-\frac{N}2}
$$
with some constant $A^*>0$, where
$$
\mathcal X=\{\varphi\in L^\infty(\mathbb R^+; H^2(\mathbb R^N)\cap L^1(\mathbb R^N)); \nabla \varphi \in L^2(\mathbb R^+; H^2(\mathbb R^N); \varphi_t\in L^2(\mathbb R^+; L^2(\mathbb R^N))\}.
$$
\end{theorem}

\section{Long-time dynamic behavior and decay rate estimation}
In this section, we are devoted to investigating the long-time asymptotic behavior of strong solutions to problem \eqref{1.6}.
 To this end, we begin by recalling following results  which come from \cite{BM} and \cite{HD}, respectively.
\begin{lemma}\label{lem2-1}(Gagliardo-Nirenberg inequality in $\mathbb R^N$)
For function $u: \mathbb R^N \to \mathbb R$, we have
$$
\|D^ju\|_{L^p(\mathbb R^N)}\le C\|D^mu\|_{L^r(\mathbb R^N)}^\alpha\|u\|^{1-\alpha}_{L^q(\mathbb R^N)}.
$$
where $1\le q, r\le \infty$, $\frac{j}m\le\alpha\le 1$,
$
\frac1p=\frac jN+(\frac1r-\frac mN)\alpha+\frac{1-\alpha}q.
$
\end{lemma}
Next, we introduce the Gagliardo-Nirenberg inequality on ball $B_R$  centered at 0 with radius $R\geq 1$. The crucial aspect is to elucidate the relationship between the embedding constants and the ball radius.

\begin{lemma}\label{GN2-2}(Gagliardo-Nirenberg inequality in a bounded domain)
For function $u: B_R\to \mathbb R$ with $B_R\subset \mathbb R^N$,
we have
$$
\|D^ju\|_{L^p(B_R)}\le C_1\|D^mu\|_{L^r(B_R)}^\alpha\|u\|^{1-\alpha}_{L^q(B_R)}+C_2R^{\frac Np-\frac Ns-j}\|u\|_{L^s(B_R)}.
$$
where $C_1>0, C_2>0$ are independent of $R$, $1\le q, r\le \infty$, $\frac{j}m\le\alpha\le 1$,
$
\frac1p=\frac jN+(\frac1r-\frac mN)\alpha+\frac{1-\alpha}q,
$
and  $s>0$ is arbitrary.
\end{lemma}

{\it\bfseries Proof.} We use scaling techniques to prove the inequality. Firstly, it is well-known that  there exist constants $ C_1>0$ and  $C_2>0$ such that  the inequality
\begin{equation}\label{GNB}
\|D^j f\|_{L^p(B_1)}\le C_1\|D^mf\|_{L^r(B_1)}^\alpha\|f\|^{1-\alpha}_{L^q(B_1)}+C_2\|f\|_{L^s(B_1)}
\end{equation}
 holds for any $f\in W^{m,r}(B_1)$.
Let
$\tilde u(y)=u(Ry)=u(x)$ for $y\in B_1$.
Then
$$
\|D^j\tilde u\|_{L^p(B_1)}=R^{j-\frac Np}\|D_ju\|_{L^p(B_R)}, \quad \|D^m\tilde u\|_{L^r(B_1)}=R^{m-\frac Nr}\|D_ju\|_{L^r(B_R)}, \quad
\|\tilde u\|_{L^q(B_1)}=R^{-\frac Nq}\|u\|_{L^q(B_R)},
$$
which in conjunction with \eqref{GNB}, results in
\begin{align*}
\|D^j u\|_{L^p(B_R)}\le C_1R^{(m-\frac Nr+\frac Nq)\alpha-\frac Nq+\frac Np-j} \|D^mu\|_{L^r(B_R)}^\alpha\|u\|^{1-\alpha}_{L^q(B_R)}+C_2R^{\frac Np-\frac Ns-j}\|u\|_{L^s(B_R)}.
\end{align*}
From $\frac1p=\frac jN+(\frac1r-\frac mN)\alpha+\frac{1-\alpha}q$, we have
$$
(m-\frac Nr+\frac Nq)\alpha-\frac Nq+\frac Np-j=0,
$$
and hence obtain that
$$
\|D^j u\|_{L^p(B_R)}\le C_1 \|D^mu\|_{L^r(B_R)}^\alpha\|u\|^{1-\alpha}_{L^q(B_R)}+C_2R^{\frac Np-\frac Ns-j}\|u\|_{L^s(B_R)}.
$$\hfill$\Box$
\begin{lemma}\label{lem2-2}
Let  $T\leq\infty$, $b>0$, $\beta>0$ and $a(t)$ be a nonnegative function locally integrable on $(0, T)$. If
$u(t)$ is nonnegative and locally integrable on $(0, T)$ with
$$
u(t)\le a(t)+b\int_0^t(t-s)^{\beta-1}u(s)ds
$$
for all $t\in (0,T)$, then
$$
u(t)\le a(t)+\theta\int_0^tE_{\beta}'(\theta(t-s))a(s)ds, \quad \text{for}\ 0\le t<T,
$$
where $\theta=(b\Gamma(\beta))^{\frac{1}{\beta}}$, $E_{\beta}(z)=\sum_{n=0}^{\infty}\frac{z^{n\beta}}{\Gamma(n\beta+1)}$,
$E_{\beta}'(z)\simeq \frac{z^{\beta-1}}{\Gamma(\beta)}$ as $z\to 0^+$.
\end{lemma}

\medskip

In order to analyze the long-time behavior of solutions, we employ a scaling transformation on the solution by letting
$$
u_{\lambda}(x,t)=\lambda^Nu(\lambda x, \lambda^2t), \quad v_{\lambda}(x,t)=v(\lambda x, \lambda^2t).
$$
It is observed that if  $(u, v)$ is the solution of \eqref{1.6}, then $(u_{\lambda}, v_{\lambda})$ is the solution to
the following problem
\begin{align}
\label{2-1}\left\{
\begin{aligned}
&\frac{\partial u}{\partial t}=\Delta u-\chi\nabla\cdot(u\nabla v) \quad  \text{in}\  Q,
\\
&\frac{\partial v}{\partial t}=\Delta v+\xi\lambda^{-N}\nabla\cdot(v\nabla u) \quad \text{in}\  Q,
\\
& u(x,0)=u_{\lambda 0}(x)=\lambda^Nu_0(\lambda x), \quad v(x,0)=v_{\lambda 0}(x)=v_0(\lambda x) \quad  \text{in}\ \mathbb R^N.
\end{aligned}\right.
\end{align}
We reformulate problem  \eqref{2-1}  into the following integral equation format:
\begin{align}
\label{2-2}\left\{
\begin{aligned}
&u_\lambda(x,t)=e^{t\Delta} u_{\lambda 0}-\chi\int_0^t\nabla\cdot\left(e^{(t-s)\Delta}(u_{\lambda}\nabla v_{\lambda})\right)ds,
\\
&v_\lambda(x,t)=e^{t\Delta} v_{\lambda 0}+\xi\lambda^{-N}\int_0^t\nabla\cdot\left(e^{(t-s)\Delta}(v_{\lambda}\nabla u_{\lambda})\right)ds,
\end{aligned}\right.
\end{align}
where
$$
(e^{t\Delta} f)(x)=\int_{\mathbb R^N}G(x-y, t)f(y)dy.
$$
As the preparation of the following analysis, we recall the following $L^p-L^q$ estimates for $e^{t\Delta}$.
\begin{lemma}\label{lem2-3}
Let $1\le q\le p\le \infty$ and $f\in L^q(\mathbb R^N)$. Then
\begin{align*}
\|e^{t\Delta} f\|_{L^p}\le (4\pi t)^{-\frac{N}2\left(\frac1q-\frac1p\right)}\|f\|_{L^q},
\\
\|\nabla e^{t\Delta} f\|_{L^p}\le C t^{-\frac{N}2\left(\frac1q-\frac1p\right)-\frac12}\|f\|_{L^q},
\end{align*}
where $C$ is a positive constant depending only on $p$ and $q$.
\end{lemma}

 For convenience, we denote by $C$ some constants in the proofs below, which may be different from each other, but do not depend on time $t$.

By using the scaling technique, we can derive  following estimates.
\begin{lemma}
\label{lem2-4}
Let $N\leq 3$,  $r>N$ and $p>1$ with $1-\frac{2}{N}<\frac1p<1+\frac1r-\frac{1}N$. Then there exists constant $C=C(N,r,p)>0$ such that  for all $ t>0$, 
we have
\begin{align}\label{2-3}
\|u(t)\|_{L^p}\le  Ct^{-\frac{N}2(1-\frac 1p)}\|u_0\|_{L^1}+Ct^{\frac12-\frac{N}{2}\left(1-\frac1p+\frac1r\right)}\|u_0\|_{L^1}
\sup_{0<s<t}\|\nabla v(s)\|_{L^r}.
\end{align}
\end{lemma}

{\it\bfseries Proof.} From the first equation of \eqref{2-2}, and applying Lemma \ref{lem2-3}, we can conclude that
\begin{align*}
\|u_{\lambda}(t)\|_{L^p}
 & \le\|e^{t\Delta} u_{\lambda 0}\|_{L^p}+\chi \int_0^t\left\|\nabla\cdot\left(e^{(t-s)\Delta}(u_{\lambda}\nabla v_{\lambda})\right)\right\|_{L^p}ds
\\
&\le(4\pi t)^{-\frac{N}2(1-\frac{1}p)}\|u_{\lambda 0}\|_{L^1}+C\chi \sup_{0<t<1}\|\nabla v_{\lambda}(t)\|_{L^r}\int_0^t(t-s)^{-\frac12-\frac{N}{2r}}\|u_{\lambda}(s)\|_{L^p}ds
\end{align*}
for $0<t\le 1$.
Thanks to Lemma \ref{lem2-2}, we get
\begin{align}\label{2-4}
&\|u_{\lambda}(t)\|_{L^p} \nonumber\\
\le & (4\pi t)^{-\frac{N}2(1-\frac{1}p)}\|u_{\lambda 0}\|_{L^1}+C\sup_{0<s< 1}\|\nabla v_{\lambda}(s)\|_{L^r}\int_0^t
(t-s)^{-\frac12-\frac{N}{2r}}s^{-\frac{N}2(1-\frac{1}p)}\|u_{\lambda 0}\|_{L^1}ds
\\
\le & (4\pi t)^{-\frac{N}2(1-\frac{1}p)}\|u_{\lambda 0}\|_{L^1}+C
\sup_{0<s< 1}\|\nabla v_{\lambda}(s)\|_{L^r}\|u_{\lambda 0}\|_{L^1}t^{\frac12-\frac{N}{2r}-\frac N2(1-\frac1p)}\int_0^1
s^{-\frac{N}2(1-\frac{1}p)}(1-s)^{-\frac12-\frac{N}{2r}}ds
\nonumber
\\
\le
& (4\pi t)^{-\frac{N}2(1-\frac{1}p)}\|u_{\lambda 0}\|_{L^1}+CB\left(1-\frac{N}2(1-\frac{1}p),\frac12-\frac{N}{2r}\right)
\sup_{0<s< 1}\|\nabla v_{\lambda}(s)\|_{L^r}\|u_{\lambda 0}\|_{L^1}t^{\frac12-\frac{N}{2r}-\frac N2(1-\frac1p)}
\nonumber
\end{align}
for all $t\leq 1$,  $p>1$ with $\frac N2(1-\frac1p)<1$, which is warranted by $\frac1p>1-\frac{2}{N}$.

Now taking $t=1$ in the above inequality, we obtain that
\begin{equation}\label{2-5}
\|u_{\lambda}(1)\|_{L^p}\le  C \|u_{\lambda 0}\|_{L^1}+C \|u_{\lambda 0}\|_{L^1}\sup_{0<s< 1}\|\nabla v_{\lambda}(s)\|_{L^r}.
\end{equation}
Noticing that
$$
\|u_{\lambda 0}\|_{L^1}=\lambda^N\int_{\mathbb R^N}u_0(\lambda x)dx=\int_{\mathbb R^N}u_0(y)dy=\|u_0\|_{L^1},
$$
$$
\|u_{\lambda}(1)\|_{L^p}=\left(\int_{\mathbb R^N}\lambda^{N p}|u(\lambda x, \lambda^2)|^p dx\right)^{\frac 1p}
=\left(\int_{\mathbb R^N}\lambda^{Np-N}|u(x, \lambda^2)|^p dx\right)^{\frac 1p}
=\lambda^{N(1-\frac1p)}\|u(\lambda^2)\|_{L^p},
$$
as well as
$$
\|\nabla v_{\lambda}(s)\|_{L^r}=\left(\int_{\mathbb R^N}\lambda^{r}|\nabla v(\lambda x, \lambda^2 s)|^r dx\right)^{\frac 1r}
=\left(\int_{\mathbb R^N}\lambda^{r-N}|\nabla v(x, \lambda^2 s)|^r dx\right)^{\frac 1r}
=\lambda^{1-\frac Nr}\|\nabla v(\lambda^2 s)\|_{L^r},
$$
 the choice of $\lambda =\sqrt t$ in \eqref{2-5} then yields
$$
t^{\frac{N}2(1-\frac 1p)}\|u(t)\|_{L^p}\le C\|u_0\|_{L^1}+C t^{\frac12-\frac{N}{2r}}\|u_0\|_{L^1}
\sup_{0<s< t}\|\nabla v(s)\|_{L^r} ,
$$
which readily implies that
$$
\|u(t)\|_{L^p}\le  Ct^{-\frac{N}2(1-\frac 1p)}\|u_0\|_{L^1}+Ct^{\frac12-\frac{N}{2}\left(1-\frac1p+\frac1r\right)}\|u_0\|_{L^1}
\sup_{0<s<t}\|\nabla v(s)\|_{L^r},
$$
and thus completes the proof. \hfill $\Box$

Similarly,   we also have
\begin{lemma}
\label{lem2-5}
Let $1\le N\le 3$, $r>N$ and $(u, v)$ be the solution of \eqref{1.6}. Then for all $p>1$ with $1-\frac{2}{N}<\frac1p<1+\frac1r-\frac{1}N$, there exists constant $C=C(N,r,p)>0$ such that  for all $ t>0$,
\begin{align}\label{2-6}
\|v(t)\|_{L^p}\le  Ct^{-\frac{N}2(1-\frac 1p)}\|v_0\|_{L^1}+Ct^{\frac12-\frac{N}{2}\left(1-\frac1p+\frac1r\right)}\|v_0\|_{L^1}
\sup_{0<s<t}\|\nabla u(s)\|_{L^r}.
\end{align}
\end{lemma}

In order to derive  the optimal decay rate estimate of $u,v$ via bootstrap argument,  
 we need the following lemmas in which  $\sup_{\frac t2<s<t}\|\nabla v(s)\|_{L^r}$ and $\sup_{\frac t2<s<t}\|\nabla u(s)\|_{L^r}$ are instead of
$\sup_{0<s<t}\|\nabla v(s)\|_{L^r}$, $\sup_{0<s<t}\|\nabla u(s)\|_{L^r}$ in \eqref{2-5} and \eqref{2-6}, respectively.
However, the approaches  adapted here seems to be inaccessible in the case $N=1$.
\begin{lemma}
\label{lem2-6} Let $N=2, 3$ and $N<r\le 4$. Suppose that $\sup_{t>0}\|\nabla^2 u(t)\|_{L^2}$  and  $\sup_{t>0}\|\nabla^2 v(t)\|_{L^2}$ are bounded, then for all $p>1$ with $\frac1p>1-\frac{2}{N}$, there exists $C=C(p,r,N)>0$ such that
\begin{align}\label{2-7}
\|u(t)\|_{L^p}\le C \|u_0\|_{L^1}t^{-\frac{N}2(1-\frac 1p)}(1+t^{-\frac{r-N}{2r}})
 +C\|u_{0}\|_{L^1}\sup_{\frac t2<s<t}\|\nabla v(s)\|_{L^r}t^{-\frac{N}2(1-\frac 1p)}(t^{\frac12-\frac N{2r}}+t^{1-\frac N{r}}),
 \\
 \label{2-8}
 \|v(t)\|_{L^p}\le C \|v_0\|_{L^1}t^{-\frac{N}2(1-\frac 1p)}(1+t^{-\frac{r-N}{2r}})
 +C\|v_{0}\|_{L^1}\sup_{\frac t2<s<t}\|\nabla u(s)\|_{L^r}t^{-\frac{N}2(1-\frac 1p)}(t^{\frac12-\frac N{2r}}+t^{1-\frac N{r}}),\
\end{align}
 for all $ t>0$.
\end{lemma}

{\it\bfseries Proof.} According to \eqref{2-4}, the following inequality
\begin{align}\label{2-9}
 &\|u_{\lambda}(t)\|_{L^q} \le (4\pi t)^{-\frac{N}2(1-\frac{1}q)}\|u_{\lambda 0}\|_{L^1}+C
\sup_{0<t< 1}\|\nabla v_{\lambda}(t)\|_{L^r}\|u_{\lambda 0}\|_{L^1}t^{\frac12-\frac{N}{2r}-\frac N2(1-\frac1 q)}
\end{align}
holds for $0<t\leq 1$, $r>N$ and $q>1$ with $\frac1 q>1-\frac{2}{N}$.

On the other hand, using Lemma \ref{lem2-3}, we have
\begin{align}
 & \|u_{\lambda}(t)\|_{L^p} \le \|e^{t\Delta} u_{\lambda 0}\|_{L^p}+\chi \int_0^t\left\|\nabla\cdot\left(e^{(t-s)\Delta}(u_{\lambda}\nabla v_{\lambda})\right)\right\|_{L^p}ds\nonumber
\\
\le &\|e^{t\Delta} u_{\lambda 0}\|_{L^p}+C\chi \int_0^{\frac t2}(t-s)^{-\frac12-\frac{N}2\left(1-\frac{1}p\right)} \|u_{\lambda}(s)\nabla v_{\lambda}(s)\|_{L^1}ds+C\chi \int_{\frac t2}^t(t-s)^{-\frac12-\frac{N}{2r}}\|u_{\lambda}(s)\|_{L^p}\|\nabla v_{\lambda}(s)\|_{L^r}ds \nonumber
\\
\le &(4\pi t)^{-\frac{N}2(1-\frac{1}p)}\|u_{\lambda 0}\|_{L^1}+C\chi\sup_{0<s<\frac t2}\|\nabla v_{\lambda}(s)\|_{L^{\frac{Nr}{2r-N}}} t^{-\frac12-\frac{N}2\left(1-\frac{1}p\right)}\int_0^{\frac t2}\|u_{\lambda}(s)\|_{L^{\frac{Nr}{Nr+N-2r}}}ds \nonumber
\\
\label{2-10}
&+C\chi\sup_{\frac t2<s<t}\|\nabla v_{\lambda}(s)\|_{L^r} \int_{\frac t2}^t(t-s)^{-\frac12-\frac{N}{2r}}\|u_{\lambda}(s)\|_{L^p}ds
\end{align}
for $0<t\leq 1$.

Now as the application of  \eqref{2-9} to $q=\frac{Nr}{Nr+N-2r}$ and $q=p$ respectively, one can have
\begin{align*}
\int_0^{\frac t2}\|u_{\lambda}(s)\|_{L^{\frac{Nr}{Nr+N-2r}}}ds \le & \int_0^{\frac t2}\left(
(4\pi s)^{\frac{N}{2r}-1}\|u_{\lambda 0}\|_{L^1}+C
\sup_{0<t< 1}\|\nabla v_{\lambda}(t)\|_{L^r}\|u_{\lambda 0}\|_{L^1}s^{-\frac12}\right)ds
\\
\le & C\|u_{\lambda 0}\|_{L^1}t^{\frac{N}{2r}}+C
\sup_{0<t< 1}\|\nabla v_{\lambda}(t)\|_{L^r}\|u_{\lambda 0}\|_{L^1}t^{\frac12}
\end{align*}
and
\allowdisplaybreaks
\begin{align*}
& \int_{\frac t2}^t(t-s)^{-\frac12-\frac{N}{2r}}\|u_{\lambda}(s)\|_{L^p}ds
\\
\le & \|u_{\lambda 0}\|_{L^1}\int_{\frac t2}^t(t-s)^{-\frac12-\frac{N}{2r}}\left((4\pi s)^{-\frac{N}2(1-\frac{1}p)} +C
\sup_{0<t<1}\|\nabla v_{\lambda}(t)\|_{L^r}s^{\frac12-\frac{N}{2r}-\frac N2(1-\frac1p)}\right)ds
\\
=&\|u_{\lambda 0}\|_{L^1}(4\pi)^{-\frac{N}2(1-\frac{1}p)}t^{\frac12-\frac{N}{2r}-\frac{N}2(1-\frac{1}p)}\int_{\frac12}^1
(1-s)^{-\frac12-\frac{N}{2r}}s^{-\frac{N}2(1-\frac{1}p)} ds
\\
&+C\|u_{\lambda 0}\|_{L^1}
\sup_{0<t<1}\|\nabla v_{\lambda}(t)\|_{L^r}t^{1-\frac{N}{r}-\frac N2(1-\frac1p)}
\int_{\frac12}^1
(1-s)^{-\frac12-\frac{N}{2r}}s^{\frac12-\frac{N}{2r}-\frac N2(1-\frac1p)}ds
\\
\le &C\|u_{\lambda 0}\|_{L^1}t^{\frac12-\frac{N}{2r}-\frac{N}2(1-\frac{1}p)}+C\|u_{\lambda 0}\|_{L^1}
\sup_{0<t< 1}\|\nabla v_{\lambda}(t)\|_{L^r}t^{1-\frac{N}{r}-\frac N2(1-\frac1p)}
\end{align*}
for $0<t\le 1$, due to $r>N$.

Substituting the above two inequalities into \eqref{2-10} yields
\begin{align*}
\|u_{\lambda}(t)\|_{L^p}
\le &  (4\pi t)^{-\frac{N}2(1-\frac{1}p)}\|u_{\lambda 0}\|_{L^1}+C\|u_{\lambda 0}\|_{L^1}\sup_{0<s<\frac t2}\|\nabla v_{\lambda}(s)\|_{L^{\frac{Nr}{2r-N}}} t^{-\frac12-\frac{N}2(1-\frac{1}p)+\frac{N}{2r}}
\\ &+C\|u_{\lambda 0}\|_{L^1}\sup_{0<s<\frac t2}\|\nabla v_{\lambda}(s)\|_{L^{\frac{Nr}{2r-N}}}
\sup_{0<s< 1}\|\nabla v_{\lambda}(s)\|_{L^r}t^{-\frac{N}2(1-\frac{1}p)} \nonumber
\\
&+C\|u_{\lambda 0}\|_{L^1} \sup_{\frac t2<s<t}\|\nabla v_{\lambda}(s)\|_{L^r}\left(t^{\frac12-\frac{N}{2r}-\frac{N}2(1-\frac{1}p)}+
\sup_{0<t< 1}\|\nabla v_{\lambda}(t)\|_{L^r}t^{1-\frac{N}{r}-\frac N2(1-\frac1p)}\right).
\end{align*}
So taking $t=1$ in the above inequality, we have
\begin{align}\label{2-11}
\|u_{\lambda}(1)\|_{L^p}\le & C \|u_{\lambda 0}\|_{L^1}+C \|u_{\lambda 0}\|_{L^1}\sup_{0<s<\frac 12}\|\nabla v_{\lambda}(s)\|_{L^{\frac{Nr}{2r-N}}} \left(1+
\sup_{0<s< 1}\|\nabla v_{\lambda}(s)\|_{L^r}\right) \nonumber
\\
&+C\|u_{\lambda 0}\|_{L^1} \sup_{\frac 12<s<1}\|\nabla v_{\lambda}(s)\|_{L^r}\left(1+
\sup_{0<s< 1}\|\nabla v_{\lambda}(s)\|_{L^r} \right).
\end{align}
Letting $\lambda =\sqrt t$ in \eqref{2-11}
and noticing that
$$
\|u_{\lambda 0}\|_{L^1}=\|u_0\|_{L^1},
\quad
\|u_{\lambda}(1)\|_{L^p}=\lambda^{N(1-\frac1p)}\|u(\lambda^2)\|_{L^p},
\quad
\|\nabla v_{\lambda}(t)\|_{L^q}
=\lambda^{1-\frac Nq}\|\nabla v(\lambda^2t)\|_{L^q},
$$
we arrive at
\begin{align*}
t^{\frac{N}2(1-\frac 1p)}\|u(t)\|_{L^p}\le & C \|u_0\|_{L^1}+C\|u_0\|_{L^1}\sup_{0<s<\frac t2}\|\nabla v(s)\|_{L^{\frac{Nr}{2r-N}}} \left(1+
t^{\frac12-\frac N{2r}}\sup_{0<s< t}\|\nabla v(s)\|_{L^r}\right)t^{\frac{N-r}{2r}} \nonumber
\\
&+C\|u_{0}\|_{L^1} t^{\frac12-\frac N{2r}}\sup_{\frac t2<s<t}\|\nabla v(s)\|_{L^r}\left(1+
t^{\frac12-\frac N{2r}}\sup_{0<s< t}\|\nabla v(s)\|_{L^r} \right),
\end{align*}
that is
\begin{align}\label{2-12}
&\|u(t)\|_{L^p}\le C \|u_0\|_{L^1}t^{-\frac{N}2(1-\frac 1p)}+C\|u_0\|_{L^1}\sup_{0<s<\frac t2}\|\nabla v(s)\|_{L^{\frac{Nr}{2r-N}}}t^{-\frac{N}2(1-\frac 1p)-\frac{r-N}{2r}}\nonumber
\\
&+C\|u_0\|_{L^1}\sup_{0<s<\frac t2}\|\nabla v(s)\|_{L^{\frac{Nr}{2r-N}}}\sup_{0<s< t}\|\nabla v(s)\|_{L^r}t^{-\frac{N}2(1-\frac 1p)}
+C\|u_{0}\|_{L^1}\sup_{\frac t2<s<t}\|\nabla v(s)\|_{L^r}t^{\frac12-\frac N{2r}-\frac{N}2(1-\frac 1p)}
\\
&+C\|u_{0}\|_{L^1} \sup_{\frac t2<s<t}\|\nabla v(s)\|_{L^r}\sup_{0<s< t}\|\nabla v(s)\|_{L^r}t^{1-\frac N{r}-\frac{N}2(1-\frac 1p)}.\nonumber
\end{align}
In addition, by Lemma \ref{lem2-1} (the Gagliardo-Nirenberg inequality),  there exists $C>0$ such that for any $q\in [\frac43, \frac{2N}{(N-2)_+})$,
\begin{equation}\label{2-13}
\|\nabla v\|_{L^q}\le C\|v\|_{L^1}^{1-\alpha}\|\nabla^2 v\|_{L^2}^{\alpha},
\end{equation}
where $\alpha=\frac{1+\frac1N-\frac1 q}{1+\frac2N-\frac12}\ge \frac12$. Noticing that
$\frac 43\leq \frac{Nr}{2r-N} < r< \frac{2N}{(N-2)_+}$ due to $N<r\le 4$, hence, $\sup_{0<s<\frac t2}\|\nabla v(s)\|_{L^{\frac{Nr}{2r-N}}}$ and $\sup_{0<s<t}\|\nabla v(s)\|_{L^{r}}$  are bounded since
$\sup_{t>0}\|\nabla^2 v(t)\|_{L^2}$ is bounded.
 Therefore combining  \eqref{2-12} with \eqref{2-13}, we readily have
\begin{align*}
\|u(t)\|_{L^p}\le C \|u_0\|_{L^1}t^{-\frac{N}2(1-\frac 1p)}(1+t^{-\frac{r-N}{2r}})
 +C\|u_{0}\|_{L^1}\sup_{\frac t2<s<t}\|\nabla v(s)\|_{L^r}t^{-\frac{N}2(1-\frac 1p)}(t^{\frac12-\frac N{2r}}+t^{1-\frac N{r}}),
\end{align*}
which is the desired inequality \eqref{2-7}.  It is observed that the inequality \eqref{2-8}  can be proved similarly.
\hfill $\Box$

 \medskip

At this position, we utilize Lemma \ref{lem2-6} to derive the optimal decay rate estimate of $u,v$ which can be stated as follows.
\begin{lemma}
\label{lem2-7}
Let $N=2, 3$, and $(u, v)$ be the solution of \eqref{1.6} such that $\sup_{t>0}\|\nabla^2 u(t)\|_{L^2}$ and $\sup_{t>0}\|\nabla^2 v(t)\|_{L^2}$ are bounded.
Then there exists $C>0$ such that for all $t\ge 1$,
\begin{align}\label{2-14}\begin{aligned}
&\|u(t)\|_{L^p}\le  Ct^{-\frac{N}2(1-\frac 1p)}\|u_0\|_{L^1}
\\
&\|v(t)\|_{L^p}\le  Ct^{-\frac{N}2(1-\frac 1p)}\|v_0\|_{L^1}
\end{aligned}\,
\end{align}
is valid for  $p\in (1,2]$.
\end{lemma}

{\it\bfseries Proof.} Recalling \eqref{2-7} and \eqref{2-8}, we see that for any $N<r\le 4$, $p>1$ with $1-\frac{2}{N}<\frac1p<1-\frac{2}{N}+\frac2r$,
\begin{align}\label{2-15}
\|u(t)\|_{L^p}\le C \|u_0\|_{L^1}t^{-\frac{N}2(1-\frac 1p)}
 +C\|u_{0}\|_{L^1}\sup_{\frac t2<s<t}\|\nabla v(s)\|_{L^r}t^{-\frac{N}2(1-\frac 1p)+1-\frac N{r}}
 \end{align}
 as well as
 \begin{align}\label{2-16}
 \|v(t)\|_{L^p}\le C \|v_0\|_{L^1}t^{-\frac{N}2(1-\frac 1p)}
 +C\|v_{0}\|_{L^1}\sup_{\frac t2<s<t}\|\nabla u(s)\|_{L^r}t^{-\frac{N}2(1-\frac 1p)+1-\frac N{r}}
\end{align}
for $t\ge \frac 12$.
From Lemma \ref{lem2-1}, we infer that
\begin{align}\label{2-17}
\|\nabla v\|_{L^r}\le C\|v\|_{L^p}^{\alpha}\|\nabla^2 v\|_{L^2}^{1-\alpha}
\end{align}
with $\alpha=\frac{\frac1N+\frac1r-\frac12}{\frac1p+\frac2N-\frac12}$. It is observed that $\alpha< \frac 12$ can be warranted by
$
p<\frac{2r}{(4-r)_+}.
$

Noticing that $\|\nabla^2 v\|_{L^2}$ is bounded, substituting \eqref{2-17} into \eqref{2-15} then yields
\begin{align}\label{2-18}
\|u(t)\|_{L^p}\le  Ct^{-\frac{N}2(1-\frac 1p)}\|u_0\|_{L^1}+Ct^{-\frac{N}2(1-\frac 1p)+1-\frac N{r}}\|u_0\|_{L^1}\sup_{\frac t2<s<t} \|v(s)\|_{L^p}^{\alpha}
\end{align}
for all $t\ge \frac 12$.
 Similarly, we  have
\begin{align}\label{2-19}
\|v(t)\|_{L^p}\le  Ct^{-\frac{N}2(1-\frac 1p)}\|v_0\|_{L^1}+Ct^{-\frac{N}2(1-\frac 1p)+1-\frac N{r}}\|v_0\|_{L^1}\sup_{\frac t2<s<t} \|u(s)\|_{L^p}^{\alpha} .
\end{align}
Furthermore,  taking $r=N+\varepsilon$, one can find $\varepsilon_1>0$ such that for $\varepsilon<\varepsilon_1$,
$1-\frac N{r}=1-\frac N{N+\varepsilon}<\frac{N}4(1-\frac 1p)$. In addition, by Lemma \ref{lem2-1},
\begin{align}\label{GN2.20}
\|u\|_{L^q}\le \|u\|_{L^1}^{1-\beta}\|D^2 u\|_{L^2}^{\beta}, \quad \text{for any}\ 1<q\le \infty,
\end{align}
with $\beta=\frac{p-1}{p(\frac2N+\frac12)}$ for all $N\le 3$. Hence from \eqref{2-19}, it follows that
\begin{align}\label{2-20}
\|v(t)\|_{L^p}\le Ct^{-\frac{N}4(1-\frac 1p)}\|v_0\|_{L^1}
\end{align}
for  $t\geq \frac12$.

Now at this position, we insert \eqref{2-20} into \eqref{2-18} to get
\begin{align}\label{2-21}
\|u(t)\|_{L^p}\le  Ct^{-\frac{N}2(1-\frac 1p)}\|u_0\|_{L^1}+Ct^{-\frac{N}2(1-\frac 1p)+1-\frac N{N+\varepsilon}
-\frac{N}4(1-\frac 1p)\cdot\frac{\frac1N+\frac1 { N+\varepsilon}-\frac12}{\frac1p+\frac2N-\frac12}
}\|u_0\|_{L^1}
\end{align}
and thus readily arrive at the decay estimate of $u$ in \eqref{2-14} for $t\ge 1$ due to the fact that
$$\displaystyle\lim_{\varepsilon\rightarrow 0 } 1-\frac N{N+\varepsilon}
-\frac{N}4(1-\frac 1p)\cdot\frac{\frac1N+\frac1 { N+\varepsilon}-\frac12}{\frac1p+\frac2N-\frac12}=
-\frac{N}4(1-\frac 1p)\cdot\frac{\frac 2N-\frac12}{\frac1p+\frac2N-\frac12}.
$$
The decay estimate of $v$ in \eqref{2-14} can be achieved similarly, and the proof of this lemma is thus completed.

\hfill $\Box$

\medskip

Our subsequent goal is to verify the result of Lemma \ref{lem2-7} for all $p$-value, inter alia  $p=\infty$. To this end, we achieve following decay result in the style of reasoning from  Lemmas  \ref{lem2-4} and \ref{lem2-6}.
\begin{lemma}
\label{lem2-8}
Assume that $N=2, 3$ and $(u, v)$ is the solution of \eqref{1.6}. If  $\sup_{t>0}\|\nabla u(t)\|_{L^r}$ and $\sup_{t>0}\|\nabla v(t)\|_{L^r}$ are bounded for some $r\in (N,4]$, then  for any  $q\in(1,r)$ with $\frac1q<\frac{2}{N}+\frac1r$, there exists a constant $C=C(q,r,N)>0$ such that
\begin{align}\label{2-23}
\|u(t)\|_{L^\infty}\le  Ct^{-\frac{N}{2q}}(1+t^{\frac{N}{2r}-\frac12})\|u_0\|_{L^q}
 +Ct^{-\frac{N}{2q}+\frac12-\frac{N}{2r}}\|u_0\|_{L^q} \sup_{\frac t2<s<t}\|\nabla v(s)\|_{L^r}\left(1+
t^{\frac12-\frac{N}{2r}}\right),
\\
\label{2-24}
\|v(t)\|_{L^\infty}\le  Ct^{-\frac{N}{2q}}(1+t^{\frac{N}{2r}-\frac12})\|v_0\|_{L^q}
 +Ct^{-\frac{N}{2q}+\frac12-\frac{N}{2r}}\|v_0\|_{L^q} \sup_{\frac t2<s<t}\|\nabla u(s)\|_{L^r}\left(1+
t^{\frac12-\frac{N}{2r}}\right)
\end{align}
for any $t>0$.
\end{lemma}

{\it\bfseries Proof.} From the first equation of \eqref{2-2}, we infer that
\begin{align*}
\|u_{\lambda}(t)\|_{L^p}\le & \|e^{t\Delta} u_{\lambda 0}\|_{L^p}+\chi \int_0^t\left\|\nabla\cdot\left(e^{(t-s)\Delta}(u_{\lambda}\nabla v_{\lambda})\right)\right\|_{L^p}ds
\\
\le &\|e^{t\Delta} u_{\lambda 0}\|_{L^p}+C\chi \int_0^t(t-s)^{-\frac12-\frac{N}{2r}}\|u_{\lambda}(s)\|_{L^p}\|\nabla v_{\lambda}(s)\|_{L^r}ds
\\
\le &(4\pi t)^{-\frac{N}{2}\left(\frac1q-\frac1p\right)}\|u_{\lambda 0}\|_{L^q}+C\chi \sup_{s\leq t}\|\nabla v_{\lambda}(s)\|_{L^r}\int_0^t(t-s)^{-\frac12-\frac{N}{2r}}\|u_{\lambda}(s)\|_{L^p}ds
\end{align*}
for $1<q<p$, $r>N$. Furthermore, with the help of Lemma \ref{lem2-2}, we can see that if $\frac{N}{2}\left(\frac1q-\frac1p\right)<1$, then for  $0<t\le 1$,
\begin{align}\label{2-25}
\|u_{\lambda}(t)\|_{L^p}
 &\le (4\pi t)^{-\frac{N}{2}\left(\frac1q-\frac1p\right)}\|u_{\lambda 0}\|_{L^q}+C\chi\sup_{0<t<1}\|\nabla v_{\lambda}(t)\|_{L^r}\int_0^t
(t-s)^{-\frac12-\frac{N}{2r}}s^{-\frac{N}{2}\left(\frac1q-\frac1p\right)}\|u_{\lambda 0}\|_{L^q}ds\nonumber
\\
&\le (4\pi t)^{-\frac{N}{2}\left(\frac1q-\frac1p\right)}\|u_{\lambda 0}\|_{L^q}+C\chi
\sup_{0<t<1}\|\nabla v_{\lambda}(t)\|_{L^r}\|u_{\lambda 0}\|_{L^q}t^{\frac12-\frac{N}{2r}-\frac{N}{2}\left(\frac1q-\frac1p\right)}\int_0^1
s^{-\frac{N}2(\frac{1}q-\frac{1}p)}(1-s)^{-\frac12-\frac{N}{2r}}ds\nonumber
\\
&\le (4\pi t)^{-\frac{N}{2}\left(\frac1q-\frac1p\right)}\|u_{\lambda 0}\|_{L^q}+C\chi
\sup_{0<t<1}\|\nabla v_{\lambda}(t)\|_{L^r}\|u_{\lambda 0}\|_{L^q}t^{\frac12-\frac{N}{2r}-\frac{N}{2}\left(\frac1q-\frac1p\right)}.
\end{align}
Similarly, the following inequality is also valid
\begin{align}\label{2-26}
\|v_{\lambda}(t)\|_{L^p}\le
(4\pi t)^{-\frac{N}{2}\left(\frac1q-\frac1p\right)}\|v_{\lambda 0}\|_{L^q}+C\chi\lambda^{-N}
\sup_{0<t<1}\|\nabla u_{\lambda}(t)\|_{L^r}\|v_{\lambda 0}\|_{L^q}t^{\frac12-\frac{N}{2r}-\frac{N}{2}\left(\frac1q-\frac1p\right)}.
\end{align}
On the other hand, we can also obtain
\begin{align}\label{2-27}
 \|u_{\lambda}(t)\|_{L^\infty} \le & \|e^{t\Delta} u_{\lambda 0}\|_{L^\infty}+\chi \int_0^t\left\|\nabla\cdot\left(e^{(t-s)\Delta}(u_{\lambda}\nabla v_{\lambda})\right)\right\|_{L^\infty}ds\nonumber
\\
\le &\|e^{t\Delta} u_{\lambda 0}\|_{L^\infty}+C\chi \int_0^{\frac t2}(t-s)^{-\frac32} \|u_{\lambda}(s)\nabla v_{\lambda}(s)\|_{L^{\frac N2}}ds+C\chi \int_{\frac t2}^t(t-s)^{-\frac12-\frac{N}{2r}}\|u_{\lambda}(s)\|_{L^\infty}\|\nabla v_{\lambda}(s)\|_{L^r}ds \nonumber
\\
\le &(4\pi t)^{-\frac{N}{2q} }\|u_{\lambda 0}\|_{L^q}+C\chi\sup_{0<s<\frac t2}\|\nabla v_{\lambda}(s)\|_{L^{\frac{Nr}{2r-N}}} t^{-\frac32}\int_0^{\frac t2}\|u_{\lambda}(s)\|_{L^r}ds \nonumber
\\
&+C\chi\sup_{\frac t2<s<t}\|\nabla v_{\lambda}(s)\|_{L^r} \int_{\frac t2}^t(t-s)^{-\frac12-\frac{N}{2r}}\|u_{\lambda}(s)\|_{L^\infty}ds,
\end{align}
where  the last two integral terms  can be estimated as follows. Indeed, by \eqref{2-25} and $\frac{N}{2}\left(\frac1q-\frac1r\right)\in (0,1)$ due to $q\in(1,r)$ with $\frac1q<\frac{2}{N}+\frac1r$, we have
\begin{align*}
\int_0^{\frac t2}\|u_{\lambda}(s)\|_{L^r}ds \le & \int_0^{\frac t2}\left((4\pi s)^{-\frac{N}{2}\left(\frac1q-\frac1r\right)}\|u_{\lambda 0}\|_{L^q}+C\chi
\sup_{0<t<1}\|\nabla v_{\lambda}(t)\|_{L^r}\|u_{\lambda 0}\|_{L^q}s^{\frac12-\frac{N}{2q}}\right)ds
\\
\le & C\|u_{\lambda 0}\|_{L^q}t^{1-\frac{N}{2}\left(\frac1q-\frac1r\right)}+C
\sup_{0<t< 1}\|\nabla v_{\lambda}(t)\|_{L^r}\|u_{\lambda 0}\|_{L^q}t^{\frac32-\frac{N}{2q}}
\end{align*}
and
\begin{align*}
\int_{\frac t2}^t(t-s)^{-\frac12-\frac{N}{2r}}\|u_{\lambda}(s)\|_{L^\infty}ds
\le  & \|u_{\lambda 0}\|_{L^q}\int_{\frac t2}^t(t-s)^{-\frac12-\frac{N}{2r}}\left((4\pi s)^{-\frac{N}{2q}} +C
\sup_{0<t<1}\|\nabla v_{\lambda}(t)\|_{L^r}s^{\frac12-\frac{N}{2r}-\frac N{2q}}\right)ds
\\
\le & C\|u_{\lambda 0}\|_{L^q} t^{\frac12-\frac{N}{2r}-\frac{N}{2q}}+C\|u_{\lambda 0}\|_{L^q}
\sup_{0<t<1}\|\nabla v_{\lambda}(t)\|_{L^r}t^{1-\frac{N}{r}-\frac N{2q}}
\end{align*}
for $0<t\le 1$. 

Substituting the above two inequalities into \eqref{2-27}, and setting $t=1$, we then get
\begin{align}
\|u_{\lambda}(1)\|_{L^\infty}
\le &  C\|u_{\lambda 0}\|_{L^q}+C\|u_{\lambda 0}\|_{L^q}\sup_{0<s<\frac 12}\|\nabla v_{\lambda}(s)\|_{L^{\frac{Nr}{2r-N}}} +C\|u_{\lambda 0}\|_{L^q}\sup_{0<s<\frac 12}\|\nabla v_{\lambda}(s)\|_{L^{\frac{Nr}{2r-N}}}
\sup_{0<s< 1}\|\nabla v_{\lambda}(s)\|_{L^r}  \nonumber
\\
\label{2-28}
&+C\|u_{\lambda 0}\|_{L^q} \sup_{\frac 12<s<1}\|\nabla v_{\lambda}(s)\|_{L^r}\left(1+
\sup_{0<t< 1}\|\nabla v_{\lambda}(t)\|_{L^r} \right).
\end{align}
Further, thanks to the fact that
$$
\|u_{\lambda 0}\|_{L^q}=\lambda^{N(1-\frac1q)}\|u_0\|_{L^q}, \quad
\|u_{\lambda}(t)\|_{L^p}=\lambda^{N(1-\frac1p)}\|u(\lambda^2t)\|_{L^p},
\quad
\|\nabla v_{\lambda}(t)\|_{L^r}
=\lambda^{1-\frac Nr}\|\nabla v(\lambda^2t)\|_{L^r},
$$
 we then choose $\lambda =\sqrt t$ in \eqref{2-28} to arrive at
\begin{align}\label{2-29}
t^{\frac N2}\|u(t)\|_{L^\infty}&\le  Ct^{\frac N2-\frac{N}{2q}}\|u_0\|_{L^q}+Ct^{\frac N2-\frac{N}{2q}+\frac{N}{2r}-\frac12}\|u_0\|_{L^q}
\sup_{0<s<\frac t2}\|\nabla v(s)\|_{L^{\frac{Nr}{2r-N}}}
\\
&+Ct^{\frac N2-\frac{N}{2q}}\|u_0\|_{L^q}\sup_{0<s<\frac t2}\|\nabla v(s)\|_{L^{\frac{Nr}{2r-N}}}
\sup_{0<s<t}\|\nabla v(s)\|_{L^r}  \nonumber
\nonumber\\
&+Ct^{\frac N2-\frac{N}{2q}+\frac12-\frac{N}{2r}}\|u_0\|_{L^q} \sup_{\frac t2<s<t}\|\nabla v(s)\|_{L^r}\left(1+
t^{\frac12-\frac{N}{2r}}\sup_{0<s< t}\|\nabla v(s)\|_{L^r} \right).\nonumber
\end{align}
Recalling \eqref{2-13}, $\sup_{0<s<\frac t2}\|\nabla v(s)\|_{L^{\frac{Nr}{2r-N}}}$ and $\sup_{0<s< t}\|\nabla v(s)\|_{L^r}$ are bounded
since $\frac43<\frac{Nr}{2r-N}<\frac{2N}{(N-2)_+}$ if $N<r\le 4$. Therefore, from \eqref{2-29}, it follows that
\begin{align*}
\|u(t)\|_{L^\infty}\le  Ct^{-\frac{N}{2q}}(1+t^{\frac{N}{2r}-\frac12})\|u_0\|_{L^q}
 +Ct^{-\frac{N}{2q}+\frac12-\frac{N}{2r}}\|u_0\|_{L^q} \sup_{\frac t2<s<t}\|\nabla v(s)\|_{L^r}\left(1+
t^{\frac12-\frac{N}{2r}}\right).
\end{align*}
The proof of \eqref{2-24} can be proved  similarly, and the proof of this lemma is thus completed.
\hfill $\Box$

\medskip
Similarly as in the proof Lemma \ref{lem2-7}, one can sharp the decay results of $(u,v)$ in Lemma \ref{lem2-8},  which is stated as
\begin{lemma}
\label{lem2-9}
Assume that $N=2, 3$.
Let $(u, v)$ be the solution of \eqref{1.6} with $\sup_{t>0}\|\nabla^2 u(t)\|_{L^2}$ and $\sup_{t>0}\|\nabla^2 v(t)\|_{L^2}$ being bounded. Then for any $t\ge 1$,
\begin{align}\label{2-30}
&\|u(t)\|_{L^\infty}\le  Ct^{-\frac{N}4}\|u_0\|_{L^2},
\\
\label{2-31}
&\|v(t)\|_{L^\infty}\le  Ct^{-\frac{N}4}\|v_0\|_{L^2},
\end{align}
where $C$ is independent of $t$.
\end{lemma}

{\it\bfseries Proof.} Recalling \eqref{2-23} and \eqref{2-24} with $q=2$, we see that for any $N<r< 4$,
\begin{align}\label{2-32}
\|u(t)\|_{L^\infty}\le  Ct^{-\frac{N}{4}}\|u_0\|_{L^2}
 +Ct^{-\frac{N}{4}+1-\frac{N}{r}}\|u_0\|_{L^2} \sup_{\frac t2<s<t}\|\nabla v(s)\|_{L^r}
 \\
 \label{2-33}
\|v(t)\|_{L^\infty}\le  Ct^{-\frac{N}{4}} \|v_0\|_{L^2}
 +Ct^{-\frac{N}{4}+1-\frac{N}{r}}\|v_0\|_{L^2} \sup_{\frac t2<s<t}\|\nabla u(s)\|_{L^r}
\end{align}
when $t\ge \frac12$.
From Lemma \ref{lem2-1}, there exists a constant $C>0$ such that
$$
\|\nabla v\|_{L^r}\le C\|v\|_{L^{\frac{2r}{4-r}}}^{\frac12}\|\nabla^2 v\|_{L^2}^{\frac12},
$$
which along with
$$
\|v\|_{L^{\frac{2r}{4-r}}}\le \|v\|_{L^1}^{\frac{4-r}{2r}}\|v\|_{L^\infty}^{\frac{3r-4}{2r}},
$$
leads to
\begin{equation} \label{2-34}
\|\nabla v\|_{L^r}\le C\|v\|_{L^1}^{\frac{4-r}{4r}}\|v\|_{L^\infty}^{\frac{3r-4}{4r}}\|\nabla^2 v\|_{L^2}^{\frac12}.
\end{equation}
Noticing that $\|\nabla^2 v\|_{L^2}$ is bounded, substituting  \eqref{2-34} into \eqref{2-32} then yields
\begin{align}\label{2-35}
\|u(t)\|_{L^\infty}\le  Ct^{-\frac{N}{4}}\|u_0\|_{L^2}
 +Ct^{-\frac{N}{4}+1-\frac{N}{r}}\|u_0\|_{L^2} \sup_{\frac t2<s<t}\|v\|_{L^\infty}^{\frac{3r-4}{4r}}\, \quad \text{for $t\ge \frac12$}.
\end{align}
Similarly, we also have
\begin{align}\label{2-36}
\|v(t)\|_{L^\infty}\le  Ct^{-\frac{N}{4}}\|v_0\|_{L^2}
 +Ct^{-\frac{N}{4}+1-\frac{N}{r}}\|v_0\|_{L^2} \sup_{\frac t2<s<t}\|u\|_{L^\infty}^{\frac{3r-4}{4r}},\, \quad \text{for $t\ge \frac12$}.
\end{align}
While recalling \eqref{GN2.20} and \eqref{2-20},  and taking
$r=N+\varepsilon$ in \eqref{2-35} and \eqref{2-36}, one can see that if  $\varepsilon$ is sufficiently small such that $1-\frac{N}{N+\varepsilon}<\frac N 8$, then
 \eqref{2-35} and \eqref{2-36} implies that
\begin{align}\label{2-37}
\|u(t)\|_{L^\infty}\le  Ct^{-\frac{N}{8}}\|u_0\|_{L^2},
 \\
 \label{2-38}
\|v(t)\|_{L^\infty}\le Ct^{-\frac{N}{8}}\|v_0\|_{L^2}
\end{align}
for any $t\ge \frac12$, respectively. Furthermore, substituting \eqref{2-37} and \eqref{2-38} into  \eqref{2-36} and \eqref{2-35} yields
\begin{align}\label{2-45}
&\|u(t)\|_{L^\infty}\le  Ct^{-\frac{N}{4}}\|u_0\|_{L^2}
 +Ct^{-\frac{N}{4}+1-\frac{N}{r}-\frac{N}{8}\frac{3r-4}{4r}}\|u_0\|_{L^2} \quad \text{for $t\ge 1$},
\\
\label{2-46}
&\|v(t)\|_{L^\infty}\le  Ct^{-\frac{N}{4}}\|v_0\|_{L^2}
 +Ct^{-\frac{N}{4}+1-\frac{N}{r}-\frac{N}{8}\frac{3r-4}{4r}}\|v_0\|_{L^2} \quad \text{for $t\ge 1$}
\end{align}
with $r=N+\varepsilon$, respectively. Noticing that
$$
\lim_{\varepsilon\to 0^+}1-\frac{N}{N+\varepsilon}-\frac{N}{8}\frac{3(N+\varepsilon)-4}{4(N+\varepsilon)}=-\frac{3N-4}{32}<0,
$$
one can see that for the suitably small $\varepsilon>0$, $-\frac{N}{4}+1-\frac{N}{r}-\frac{N}{8}\frac{3r-4}{4r}<-\frac{N}{4}$ with $r=N+\varepsilon$ and thereby achieve the desired inequalities \eqref{2-30} and \eqref{2-31}. \hfill $\Box$

\medskip

Now with the help of Lemma \ref{lem2-7} and Lemma \ref{lem2-9}, one  can obtain the optimal decay rate of $u$ and $v$.
\begin{lemma}
\label{lem2-10} Let the assumption of Lemma \ref{lem2-9} hold. Then
there exists $C>0$ such that
\begin{align}\label{2-41}
&\|u(t)\|_{L^p}\le  C(t+1)^{-\frac{N}2\left(1-\frac1p\right)}\|u_0\|_{L^1},
\\
&\|v(t)\|_{L^p}\le  C(t+1)^{-\frac{N}2\left(1-\frac1p\right)}\|v_0\|_{L^1}
\end{align}
for all $p\in (1,\infty]$ and $t>0$.
\end{lemma}

{\it\bfseries Proof.} From Lemma \ref{lem2-9},  it is observed that for any $t>2$,
\begin{align*}
&\|u(t)\|_{L^\infty}\le  Ct^{-\frac{N}4}\left\|u\left(\frac t2\right)\right\|_{L^2},\quad \|v(t)\|_{L^\infty}\le  Ct^{-\frac{N}4}\left\|v\left(\frac t2\right)\right\|_{L^2}.
\end{align*}
On the other hand, by leveraging Lemma \ref{lem2-7}, one can conclude that for any  $t> 2$,
\begin{align*}
&\left\|u\left(\frac t2\right)\right\|_{L^2}\le  Ct^{-\frac{N}4}\|u_0\|_{L^1},
\left\|v\left(\frac t2\right)\right\|_{L^2}\le  Ct^{-\frac{N}4}\|v_0\|_{L^1}.
\end{align*}
Hence as the combination of the above inequalities, one deduces that for any $t > 2$,
\begin{align}\label{2-43}\begin{aligned}
&\|u(t)\|_{L^\infty}\le  Ct^{-\frac{N}2}\|u_0\|_{L^1},
\\
&\|v(t)\|_{L^\infty}\le  Ct^{-\frac{N}2}\|v_0\|_{L^1},
\end{aligned}\end{align}
which implies that
\begin{align}\label{2-44}\begin{aligned}
&\|u(t)\|_{L^\infty}\le  C(t+1)^{-\frac{N}2}\|u_0\|_{L^1},
\\
&\|v(t)\|_{L^\infty}\le  C(t+1)^{-\frac{N}2}\|v_0\|_{L^1}
\end{aligned}\end{align}
for any $t>0$, due to the fact that  both $\|u(t)\|_{L^\infty}$ and $\|u(t)\|_{L^\infty}$ are bounded by Gagliardo-Nirenberg  inequality.

As a directed application of $
\|f\|_{L^p}\le \|f\|_{L^\infty}^{\frac{p-1}p} \|f\|_{L^1}^{\frac{1}p},
$
 \eqref{2-44} gives
\begin{align*}
\|u(t)\|_{L^p}\le  C(t+1)^{-\frac{N}2\left(1-\frac1p\right)}\|u_0\|_{L^1},\quad
\|v(t)\|_{L^p}\le  C(t+1)^{-\frac{N}2\left(1-\frac1p\right)}\|v_0\|_{L^1}
\end{align*}
for any $t>0$ and thus completes the proof. \hfill $\Box$

\vskip2mm
 In contrast to the case of $N=2,3$,  the decay estimate in  Lemma \ref{lem2-6} seems to be invalid for the one-dimensional case,  and thus the optimal decay rate of solution  $(u,v)$ of \eqref{1.6} is likely inaccessible. 

\begin{lemma}
\label{lem2-11}
Let $N=1$ and $(u, v)$ be the solution of \eqref{1.6} with $\sup_{t>0}\|\nabla^2 u(t)\|_{L^2}$ and $\sup_{t>0}\|\nabla^2 v(t)\|_{L^2}$ being bounded. Then there exists $C>0$ such that
\begin{align}
\begin{aligned}
&\|u(t)\|_{L^p}\le  C(t+1)^{-\frac{3}8\left(1-\frac1p\right)}\|u_0\|_{L^1},
\\
&\|v(t)\|_{L^p}\le  C(t+1)^{-\frac{3}8\left(1-\frac1p\right)}\|v_0\|_{L^1}
\end{aligned}\qquad
\end{align}
for all $p\in (1,\infty]$ and $t>0$.
\end{lemma}

{\it\bfseries Proof.} As  the application of \eqref{2-3} to $N=1$, we can see that  for $r>1$ there exists $C>0$  such that
\begin{align}
\|u(t)\|_{L^\infty}\le  Ct^{-\frac{1}2}\|u_0\|_{L^1}+Ct^{-\frac1{2r}}\|u_0\|_{L^1}
\sup_{0<s<t}\|\nabla v(s)\|_{L^r}.
\end{align}
In addition, according to Lemma \ref{lem2-1}, the inequality
$
\|\nabla v\|_{L^{r}}\le \|v\|_{L^1}^{1-\alpha}\|\nabla^2 v\|_{L^2}^{\alpha}
$
with $\alpha=\frac45-\frac2{5r}$ is valid for all $r\ge \frac43$. 

 Taking $r=\frac43$ in \eqref{2-46}, we arrive at
$$
\|u(t)\|_{L^\infty}\le  Ct^{-\frac{1}2}\|u_0\|_{L^1}+Ct^{-\frac3{8}}\|u_0\|_{L^1}
$$
and thus  for any $t>0$,
$$
\|u(t)\|_{L^\infty}\le C(t+1)^{-\frac3{8}}\|u_0\|_{L^1}
$$
since $\|u(t)\|_{L^\infty}$ is bounded, which readily implies that
$$
\|u\|_{L^p} \le C(t+1)^{-\frac3{8}(1-\frac1p)}\|u_0\|_{L^1}
$$
since
$
\|u\|_{L^p}\le \|u\|_{L^\infty}^{\frac{p-1}p}\|u\|_{L^1}^{\frac 1p}.
$
Similar to the above proof, and using Lemma \ref{lem2-5}, we also have
$$
\|v\|_{L^p} \le C(t+1)^{-\frac3{8}(1-\frac1p)}\|v_0\|_{L^1}.$$ \hfill $\Box$

In what follows we show that the solution  converges to heat kernel as time tends to infinity. As a prerequisite,  it is observed that  the  difference  between $ t^{\frac N2\left(1-\frac1p\right)} e^{t\Delta}f$ and $  t^{\frac N2\left(1-\frac1p\right)}\int_{\mathbb R^N}f(y)dy G(\cdot,t)$ decays to zero in  $L^p(\mathbb R^N)$. The proof thereof  may be  for convenient   of some readers.

\begin{lemma}
\label{lem2-12}
Assume that $f\ge 0$, and $f\in L^1(\mathbb R^N)\cap L^\infty(\mathbb R^N)$.  Then for any $1\le p\le \infty$,
\begin{align}\label{2-47}
\lim_{t\to\infty}t^{\frac N2\left(1-\frac1p\right)}\left\|e^{t\Delta}f-\int_{\mathbb R^N}f(y)dy G(\cdot,t)\right\|_{L^p}=0,
\end{align}
where $G(x,t)$  is the heat kernel.
\end{lemma}

{\it\bfseries Proof.} Given that $f(x)\in L^1(\mathbb R^N)$, then for any $\varepsilon>0$, there exists $R>0$ such that
$$
\int_{\mathbb R^N\setminus B_R}f(y)dy<\varepsilon.
$$
We also note that
$$
|G(x-y,t)-G(x,t)|=|y\cdot\nabla G(x-\theta y, t)|=(4\pi t)^{-\frac N2}\frac{|y|}{\sqrt t}\left|\frac{x-\theta y}{2\sqrt t}e^{-\frac{|x-\theta y|^2}{4t}}\right|\le \frac{1}{\sqrt 2}e^{-\frac12}(4\pi t)^{-\frac N2}\frac{|y|}{\sqrt t}.
$$
Then by a direct calculation, we see that
\begin{align*}
&(4\pi t)^{\frac N2}\left|e^{t\Delta}f-\int_{\mathbb R^N}f(y)dy G(x,t)\right|=(4\pi t)^{\frac N2}\left|\int_{\mathbb R^N}(G(x-y,t)-G(x,t))f(y)dy\right|
\\
&\le (4\pi t)^{\frac N2}\int_{\mathbb R^N\setminus B_R}|G(x-y,t)-G(x,t)|f(y)dy+(4\pi t)^{\frac N2}\int_{B_R}|G(x-y,t)-G(x,t)|f(y)dy
\\
&\le 2\int_{\mathbb R^N\setminus B_R}f(y)dy+e^{-\frac12}\frac{R}{\sqrt{2t}}\int_{B_R}f(y)dy
\\
&\le 2\varepsilon+e^{-\frac12}\frac{R}{\sqrt{2t}}\int_{\mathbb R^N}f(y)dy,
\end{align*}
which leads to
\begin{align}\label{2-48}
\lim_{t\to\infty}t^{\frac N2}\left\|e^{t\Delta}f-\int_{\mathbb R^N}f(y)dy G(\cdot,t)\right\|_{L^\infty}=0.
\end{align}

Next, we show that $\left\|e^{t\Delta}f-\int_{\mathbb R^N}f(y)dy G(\cdot,t)\right\|_{L^1}\to 0$. For this purpose, we first prove the following
inequality
\begin{align*}
\int_{\mathbb R^N}|G(x-y,t)-G(x,t)|dx& =\int_{\mathbb R^N}\left|\int_0^1\frac{d}{ds}G(x-sy,t)ds\right|dx
\\
&\le \int_{\mathbb R^N}\int_0^1\left| y\nabla G(x-sy,t)\right|dsdx
\\
&\le \int_0^1 \int_{\mathbb R^N}(4\pi t)^{-\frac N2}\frac{|y|}{\sqrt t}\left|\frac{x-s y}{2\sqrt t}e^{-\frac{|x-s y|^2}{4t}}\right|dxds
\\
&= (\pi)^{-\frac N2}\frac{|y|}{\sqrt t} \int_{\mathbb R^N}|xe^{-x^2}|dx
\\
&=(\pi)^{-\frac N2}\frac{|y|}{\sqrt t} |\partial B_1|\int_{0}^\infty r^Ne^{-r^2}dr.
\end{align*}
It follows from the latter that
\begin{align*}
&\int_{\mathbb R^N}\left|e^{t\Delta}f-\int_{\mathbb R^N}f(y)dy G(x,t)\right|dx
\\
&\le\int_{\mathbb R^N}\int_{\mathbb R^N}|G(x-y,t)-G(x,t)|f(y)dydx
\\
&\le \int_{\mathbb R^N}f(y)dy\int_{\mathbb R^N}|G(x-y,t)-G(x,t)|dx
\\
&\le C\int_{B_R}\frac{|y|}{\sqrt t}f(y)dy+ \int_{\mathbb R^N\setminus B_R}f(y)dy\int_{\mathbb R^N}(|G(x-y,t)|+|G(x,t)|)dx
\\
&\le C\frac{R}{\sqrt t}\int_{\mathbb R^N}f(y)dy+2\int_{\mathbb R^N\setminus B_R}f(y)dy
\\
&\le C\frac{R}{\sqrt{t}}\int_{\mathbb R^N}f(y)dy+2\varepsilon.
\end{align*}
As a result of this, we have
\begin{align}\label{2-49}
\lim_{t\to\infty} \left\|e^{t\Delta}f-\int_{\mathbb R^N}f(y)dy G(\cdot,t)\right\|_{L^1}=0.
\end{align}
Therefore according to \eqref{2-48} and \eqref{2-49},  \eqref{2-47} can be identified by the fact that
$$
\left\|e^{t\Delta}f-\int_{\mathbb R^N}f(y)dy G(\cdot,t)\right\|_{L^p}\le \left\|e^{t\Delta}f-\int_{\mathbb R^N}f(y)dy G(\cdot,t)\right\|_{L^1}^{\frac1p}\left\|e^{t\Delta}f-\int_{\mathbb R^N}f(y)dy G(\cdot,t)\right\|_{L^\infty}^{1-\frac1p}.
$$
\hfill $\Box$

\medskip
\begin{lemma}
\label{lem2-13}
Assume that $N=2, 3$, and  $(u, v)$ is the solution of \eqref{1.6} with $\sup_{t>0}\|\nabla^2 u(t)\|_{L^2}$ and $\sup_{t>0}\|\nabla^2 v(t)\|_{L^2}$ being bounded.
Then
\begin{align}\label{2-50}\begin{aligned}
&\lim_{t\to\infty}t^{\frac N{2}\left(1-\frac1p\right)}\left\|u(\cdot,t)-\int_{\mathbb R^N}u_0(y)dy G(\cdot, t)\right\|_{L^p}=0,
\\
&\lim_{t\to\infty}t^{\frac N{2}\left(1-\frac1p\right)}\left\|v(\cdot,t)-\int_{\mathbb R^N}v_0(y)dy G(\cdot, t)\right\|_{L^p}=0,
\end{aligned}
\end{align}
where $p\in [1,\infty]$.
\end{lemma}

{\it\bfseries Proof.} According to
$$
u(x,t)=e^{t\Delta} u_0-\chi\int_0^t\nabla\cdot\left(e^{(t-s)\Delta}(u \nabla v)\right)ds,
$$
for any $t>2$, using \eqref{2-13} and \eqref{2-41}
\begin{align*}
&\|u(\cdot,t)-e^{t\Delta} u_0\|_{L^\infty}
\\
&\le \chi\int_0^t\left\|\nabla\cdot\left(e^{(t-s)\Delta}(u \nabla v)\right)\right\|_{L^\infty}ds
\\
&\le C\chi\int_0^{\frac t2}(t-s)^{-\frac12-\frac N{2}}\|u(s)\|_{L^4}\|\nabla v(s)\|_{L^{\frac43}}ds+
 C\chi\int_{\frac t2}^t(t-s)^{-\frac12-\frac N{2r}}\|u(s)\|_{L^\infty}\|\nabla v(s)\|_{L^r}ds
 \\
 &\le C\chi\int_0^{\frac t2}(t-s)^{-\frac12-\frac N{2}}(1+s)^{-\frac{3N}8}\|u_0\|_{L^1}ds+
 C\chi\int_{\frac t2}^t(t-s)^{-\frac12-\frac N{2r}}s^{-\frac N2}\|\nabla v(s)\|_{L^r}\|u_0\|_{L^1}ds
 \\
 &\le Ct^{-\frac12-\frac N{2}}\left((2+t)^{1-\frac{3N}8}+1\right)+C \int_{\frac t2}^t(t-s)^{-\frac12-\frac N{2r}}s^{-\frac N2}\|\nabla v(s)\|_{L^r}ds.
\end{align*}
By Lemma \ref{lem2-1}, for $r\in [4, \frac{2N}{N-2}]$,
$$
\|\nabla v\|_{L^r}\le C\|v\|_{L^\infty}^{\alpha}\|D^2v\|_{L^2}^{1-\alpha}
$$
with $\alpha=\frac{\frac1N+\frac1r-\frac12}{\frac2N-\frac12}\le \frac12$, which along with Lemma \ref{lem2-10}, implies that
$$
\|\nabla v(s)\|_{L^r}\le C\|v\|_{L^\infty}^{\alpha}\le Cs^{-\frac N2 \alpha }
$$
for $s>1$.
Hence one can conclude that  for any $t>2$,
\begin{align*}
\|u(\cdot,t)-e^{t\Delta} u_0\|_{L^\infty}\le Ct^{-\frac12-\frac N{2}}\left((2+t)^{1-\frac{3N}8}+1\right)+C\int_{\frac t2}^t(t-s)^{-\frac12-\frac N{2r}}s^{-\frac N2-\frac N2 \alpha }ds
\\
\le Ct^{-\frac12-\frac N{2}}\left((2+t)^{1-\frac{3N}8}+1\right)+Ct^{\frac12-\frac N{2r}-\frac N2-\frac N2 \alpha}\int_{\frac12}^1(1-s)^{-\frac12-\frac N{2r}}s^{-\frac N2-\frac N2 \alpha }ds,
\end{align*}
which means that
\begin{align}\label{2-51}
t^{\frac N{2}}\|u(\cdot,t)-e^{t\Delta} u_0\|_{L^\infty}\le Ct^{-\frac12}\left((2+t)^{1-\frac{3N}8}+1\right)+Ct^{\frac12-\frac N{2r}-\frac N2 \alpha }.
\end{align}
Furthermore, when fixing $r=5$, we have $\alpha=\frac25$ for $N=2$ and $\alpha=\frac15$ for $N=3$.
The latter  ensures that  $\frac12-\frac N{2r}-\frac N2 \alpha=-\frac1{10}$ for $N=2, 3$.
Therefore from \eqref{2-51},  it thereby follows that
\begin{equation}\label{2-52}
t^{\frac N{2}}\|u(\cdot,t)-e^{t\Delta} u_0\|_{L^\infty}\le 
Ct^{-\frac1{10}}, \quad\text{for $t>2$}.
\end{equation}
On the other hand, similar to the above proof, we use \eqref{2-13} and \eqref{2-41} again to see that
\begin{align}
&\|u(\cdot,t)-e^{t\Delta} u_0\|_{L^1}\le \chi\int_0^t\left\|\nabla\cdot\left(e^{(t-s)\Delta}(u \nabla v)\right)\right\|_{L^1}ds\nonumber
\\
&\le C\chi\int_0^{\frac t2}(t-s)^{-\frac12}\|u(s)\|_{L^4}\|\nabla v(s)\|_{L^{\frac 43}}ds+
 C\chi\int_{\frac t2}^t(t-s)^{-\frac12}\|u(s)\|_{L^4}\|\nabla v(s)\|_{L^{\frac 43}}ds\nonumber
 \\
 &\le Ct^{-\frac12}\left((2+t)^{1-\frac{3N}8}+1\right)+C\int_{\frac t2}^t(t-s)^{-\frac12}s^{-\frac{3N}8}ds\nonumber
 \\
 \label{2-63}
 &\le Ct^{-\frac12}\left((2+t)^{1-\frac{3N}8}+1\right)+Ct^{-\frac{3N-4}8}
\end{align}
for any $t>2$,
which together with \eqref{2-52} yields
\begin{equation}\label{2-54}
\lim_{t\to\infty}t^{\frac N{2}\left(1-\frac1p\right)}\|u(\cdot,t)-e^{t\Delta} u_0\|_{L^p}=0
\end{equation}
for any $1\le p\le \infty$ since $\|f\|_{L^p}\le \|f\|_{L^\infty}^{1-\frac1p}\|f\|_{L^1}^{\frac1p}$.
Combining \eqref{2-54} and \eqref{2-47}, we can conclude that
\begin{equation}\label{2-65}
\lim_{t\to\infty}t^{\frac N{2}\left(1-\frac1p\right)}\left\|u(\cdot,t)-\int_{\mathbb R^N}u_0(y)dy G(\cdot, t)\right\|_{L^p}=0.
\end{equation}
Similarly,
\begin{equation}\label{2-66}
\lim_{t\to\infty}t^{\frac N{2}\left(1-\frac1p\right)}\left\|v(\cdot,t)-\int_{\mathbb R^N}v_0(y)dy G(\cdot, t)\right\|_{L^p}=0
\end{equation}
for $1\le p\le \infty$. This lemma is proved. \hfill $\Box$

\medskip

Upon the decay results for the one-dimensional space asserted in  Lemma \ref{lem2-11},  we have
\begin{lemma}
\label{lem2-14}
Let $N=1$, and $(u, v)$ be the strong solution of \eqref{1.6} such that  $\sup_{t>0}\|\nabla^2 u(t)\|_{L^2}$ and $\sup_{t>0}\|\nabla^2 v(t)\|_{L^2}$ are bounded.
Then for any $0<m<\frac{9}{32}$, we have
\begin{align}\label{2-67}
\lim_{t\to\infty}t^{m\left(1-\frac1p\right)}\left\|u(\cdot,t)-\int_{\mathbb R^N}u_0(y)dy G(\cdot, t)\right\|_{L^p}=0,
\end{align}
where $p\in (1,\infty]$.
\end{lemma}

{\it\bfseries Proof.} Similar to the proof of Lemma \ref{lem2-13}, using \eqref{2-13} and \eqref{2-45}, we arrive at
\begin{align*}
&\|u(\cdot,t)-e^{t\Delta} u_0\|_{L^\infty}\le \chi\int_0^t\left\|\nabla\cdot\left(e^{(t-s)\Delta}(u \nabla v)\right)\right\|_{L^\infty}ds
\\
&\le C\chi\int_0^{\frac t2}(t-s)^{-1}\|u(s)\|_{L^4}\|\nabla v(s)\|_{L^{\frac43}}ds+
 C\chi\int_{\frac t2}^t(t-s)^{-\frac12-\frac12\left(\frac1p+\frac{p+2}{4p}\right)}\|u(s)\|_{L^p}\|\nabla v(s)\|_{L^{\frac{4p}{p+2}}}ds
 \\
 &\le C\sup_{s>0}\|\nabla v(s)\|_{L^{\frac43}}\int_0^{\frac t2}(t-s)^{-1}(1+s)^{-\frac{9}{32}}ds+C\chi\int_{\frac t2}^t(t-s)^{-\frac12-\frac12\left(\frac1p+\frac{p+2}{4p}\right)}\|u(s)\|_{L^p}\|v(s)\|_{L^p}^{\frac12}\|D^2 v(s)\|_{L^2}^{\frac12}ds
 \\
 &\le C(1+t)^{-\frac{9}{32}}+C\int_{\frac t2}^t(t-s)^{-\frac{5p+6}{8p}}s^{-\frac9{16}\left(1-\frac1p\right)}ds
 \\
 &\le  Ct^{-\frac{9}{32}}+Ct^{1-\frac{5p+6}{8p}-\frac9{16}\left(1-\frac1p\right)}
 \\
 &\le  Ct^{-\frac{9}{32}}+Ct^{-\frac{3p+3}{16p}},
\end{align*}
for $t>2$ and  any $p>2$. It is easy to see that $\frac{3p+3}{16p}$ is decreasing on $p$ and $\frac{3p+3}{16p}\to \frac{9}{32}$ as  $p\searrow 2$. Hence for any  $0<m<\frac{9}{32}$, we have
\begin{equation}
\label{2-68}
\lim_{t\to\infty}t^{m}\|u(\cdot,t)-e^{t\Delta} u_0\|_{L^\infty}=0.
\end{equation}
Combining \eqref{2-47} and  \eqref{2-68}, we arrive at
\begin{equation}
\label{2-69}
\lim_{t\to\infty}t^{m}\left\|u(\cdot,t)-\int_{\mathbb R}u_0(y)dy G(\cdot, t)\right\|_{L^\infty}=0,
\end{equation}
and hence \eqref{2-67} results from the inequality
\begin{align*}
&\left\|u(\cdot,t)-\int_{\mathbb R}u_0(y)dy G(\cdot, t)\right\|_{L^p}
\\
&\le \left\|u(\cdot,t)-\int_{\mathbb R}u_0(y)dy G(\cdot, t)\right\|_{L^1}^{\frac1p}
\left\|u(\cdot,t)-\int_{\mathbb R}u_0(y)dy G(\cdot, t)\right\|_{L^\infty}^{1-\frac1p}
\\
&\le 2^{\frac1p}\|u_0\|_{L^1}^{\frac1p}\left\|u(\cdot,t)-\int_{\mathbb R}u_0(y)dy G(\cdot, t)\right\|_{L^\infty}^{1-\frac1p}.
\end{align*}
\hfill $\Box$

\medskip

In conclusion, Theorem \ref{thm-1} and Theorem \ref{thm-2} are the consequence of  combining  Lemma \ref{lem2-10} with Lemma \ref{lem2-13},  and  Lemma \ref{lem2-13} with Lemma \ref{lem2-14}, respectively.

\section{Existence of global strong solutions}

In Section 2, we have investigated the long-time dynamic behavior of the global strong solution $(u, v)\in L^\infty(\mathbb R^+; H^2(\mathbb R^N)\cap L^1(\mathbb R^N))$ of \eqref{1.6}. Hence this section  concentrates on  the existence of  such strong solutions to the Cauchy problem of \eqref{1.6}. For this purpose, we consider the following approximation problem of \eqref{1.6}
\begin{align}
\label{3-1}\left\{
\begin{aligned}
&u_t+\varepsilon\Delta^2u-\Delta u=-\chi\nabla\cdot(u\nabla v) && \text{in}\  Q,
\\
&v_t-\Delta v=\xi\nabla\cdot(v\nabla u) && \text{in}\  Q,
\\
& u(x,0)=u_{\varepsilon 0}(x)\ge 0,  v(x,0)=v_{\varepsilon 0}(x)\ge 0, &&  x\in\mathbb R^N,
\end{aligned}\right.
\end{align}
where $Q=\mathbb R^N\times \mathbb R^+$, $\varepsilon\in(0,1)$, $u_{\varepsilon 0}, v_{\varepsilon 0} \in C_0^\infty(\mathbb R^N)$ are the smooth approximations of $u_0$ and $v_0$ in $H^2(\mathbb R^N)$, respectively.
Let
$$
\mathcal D=L^\infty((0, T); H^2(\mathbb R^N)\cap L^1(\mathbb R^N))\cap L^2((0, T); H^3(\mathbb R^N)),
$$
$$
\tilde{\mathcal D}=\{u\in L^\infty((0, T); H^3(\mathbb R^N)\cap L^1(\mathbb R^N))\cap L^2((0, T); H^5(\mathbb R^N)),  u_t\in L^2(Q_T), xu\in L^\infty((0, T); L^2(\mathbb R^N))\}.
$$
Define
\begin{align*}
\mathcal F:  \mathcal D \times [0,1]\to \mathcal D,  \qquad \mathcal F(\tilde u, \sigma)=u,
\end{align*}
where $u$ is the solution of
\begin{align}
\label{3-2}\left\{
\begin{aligned}
&u_t+\varepsilon\Delta^2u-\Delta u=-\chi\nabla\cdot(u\nabla v) && \text{in}\  Q,
\\
&v_t-\Delta v=\xi\nabla\cdot(v\nabla\tilde u) && \text{in}\  Q,
\\
& u(x,0)=\sigma u_{\varepsilon 0}(x)\ge 0,  v(x,0)=\sigma v_{\varepsilon 0}(x)\ge 0, &&  x\in\mathbb R^N.
\end{aligned}\right.
\end{align}
According to  Proposition \ref{appro3-1}, Proposition \ref{appro3-2} and Lemma \ref{aplem3-3} in Appendix, we see that for $\tilde u\in \mathcal D$, the problem \eqref{3-2} admits a unique solution $u\in \tilde{\mathcal D}$. Therefore, $\mathcal F$ is well-defined.
By Aubin-Lions Lemma and the compact results in \cite{ZS} one also see that $\mathcal F$ is a compact mapping.
Next, it is easy to see that
$$
\mathcal F(\tilde u, 0)=0.
$$
According to the Leray-Schauder fixed point theorem,
to show the existence of strong solutions to the problem \eqref{3-1}, it remains only to show that there exists constant $C>0$ independent of $\sigma\in [0,1]$ such that
$
\|u\|_{\mathcal D}\le C $
for any  $\mathcal F(u, \sigma)=u$.

\begin{lemma}\label{lem3-1}
Assume that $\mathcal F(u, \sigma)=u$.
Then we have
\begin{align}
&\frac{d}{dt} \left(\|u(\cdot, t)\|_{H^2}^{2}+\|v(\cdot, t)\|_{H^2}^{2}\right)
+2\varepsilon\|\Delta u(\cdot,t)\|_{H^2}^2+ \|\nabla u(\cdot,t)\|_{H^2}^2+\|\nabla v(\cdot,t)\|_{H^2}^2  \nonumber
\\
\label{3-3}
&\le  C^*(\chi^2 +\xi^2)\left(\|u(\cdot, t)\|_{H^2}^2+\| v(\cdot, t)\|_{H^2}^2\right)\left(\|\nabla u(\cdot, t)\|_{H^2}^2
+\|\nabla v(\cdot, t)\|_{H^2}^2\right),
\end{align}
where constant $C^*>1$ depends only on the spatial dimension $N$.
\end{lemma}

{\it\bfseries Proof.} Multiplying the first equation of \eqref{3-2} by $u$,  integrating the results over $\mathbb R^N$ yields
\begin{align*}
&\frac{1}{2}\frac{d}{dt}\int_{\mathbb R^N}u^{2}dx+\varepsilon\int_{\mathbb R^N} |\Delta u|^2dx+\int_{\mathbb R^N} |\nabla u|^2dx
\\
&=\chi\int_{\mathbb R^N} u\nabla v\nabla udx
\\
&\le \frac 12\int_{\mathbb R^N} |\nabla u|^2dx+\frac{\chi^2}2\int_{\mathbb R^N} u^2|\nabla v|^2dx.
\end{align*}
Similarly  testing the second equation of \eqref{3-2} by $v$ and integrating the  corresponding result over $\mathbb R^N$ yield
\begin{align*}
&\frac{1}{2}\frac{d}{dt}\int_{\mathbb R^N}|v|^{2}dx+\int_{\mathbb R^N} |\nabla v|^2dx
\\
&=-\xi\int_{\mathbb R^N} v\nabla v\nabla udx
\\
&\le \frac12\int_{\mathbb R^N} |\nabla v|^2dx+\frac{\xi^2}2\int_{\mathbb R^N} v^2|\nabla u|^2dx.
\end{align*}
Combining the above two equalities, we have
\begin{align}
&\frac{d}{dt}\int_{\mathbb R^N} \left(|u|^{2}+|v|^{2}
\right)dx+2\varepsilon \int_{\mathbb R^N} |\Delta u|^2dx+\int_{\mathbb R^N}\left(|\nabla u|^2+ |\nabla v|^2\right)dx\nonumber
\\
&\le \chi^2\|u\|_{L^\infty}^2\|\nabla v\|_{L^2}^2+ \xi^2\|v\|_{L^\infty}^2\|\nabla u\|_{L^2}^2\nonumber
\\
\label{3-4}
&\le C_1\chi^2\|u\|_{H^2}^2\|\nabla v\|_{L^2}^2+ C_1\xi^2\|v\|_{H^2}^2\|\nabla u\|_{L^2}^2.
\end{align}
Furthermore, multiplying the first equation of \eqref{3-2} by $-\Delta u$, and integrating the resulting equation over $\mathbb R^N$ yields
\begin{align}
&\frac{1}{2}\frac{d}{dt}\int_{\mathbb R^N}|\nabla u|^{2}dx+\varepsilon\int_{\mathbb R^N} |\nabla\Delta u|^2dx+\int_{\mathbb R^N} |\Delta u|^2dx\nonumber
\\
&=\chi\int_{\mathbb R^N} \nabla\cdot(u\nabla v)\Delta udx \nonumber
\\
\label{3-5}
&\le \frac12\int_{\mathbb R^N} |\Delta u|^2dx+\chi^2 \int_{\mathbb R^N}( |\nabla u|^2|\nabla v|^2+ u^2|\Delta v|^2)dx.
\end{align}
Similarly,  we get 
\begin{align}
&\frac{1}{2}\frac{d}{dt}\int_{\mathbb R^N}|\nabla v|^{2}dx+\int_{\mathbb R^N} |\Delta v|^2dx\nonumber
\\
&=-\xi\int_{\mathbb R^N} \nabla\cdot(v\nabla u)\Delta vdx\nonumber
\\
\label{3-6}
&\le \frac12\int_{\mathbb R^N} |\Delta v|^2dx+\xi^2 \int_{\mathbb R^N}( |\nabla u|^2|\nabla v|^2+ v^2|\Delta u|^2)dx.
\end{align}
Combining \eqref{3-5} and \eqref{3-6}, we arrive at
\begin{align}
&\frac{d}{dt}\int_{\mathbb R^N}\left( |\nabla u|^{2}+ |\nabla v|^{2}\right)dx+2\varepsilon \int_{\mathbb R^N} |\nabla\Delta u|^2dx+ \int_{\mathbb R^N} (|\Delta u|^2+|\Delta v|^2)dx \nonumber
\\
&\le 2\chi^2 \int_{\mathbb R^N}( |\nabla u|^2|\nabla v|^2+ u^2|\Delta v|^2)dx+2\xi^2 \int_{\mathbb R^N}( |\nabla u|^2|\nabla v|^2+ v^2|\Delta u|^2)dx
\nonumber
\\
&\le 2\chi^2(\|\nabla u\|_{L^4}^2\|\nabla v\|_{L^4}^2+\|u\|_{L^\infty}^2\|\Delta v\|_{L^2}^2)+
2\xi^2(\|\nabla u\|_{L^4}^2\|\nabla v\|_{L^4}^2+\|v\|_{L^\infty}^2\|\Delta u\|_{L^2}^2) \nonumber
\\
\label{3-7}
&\le C(\chi^2+\xi^2)(\|\nabla u\|_{H^1}^2+\|\nabla v\|_{H^1}^2)(\|u\|_{H^2}^2+\|v\|_{H^2}^2).
\end{align}
Moreover, applying $\nabla$ to the  equations in \eqref{3-2}, and multiplying the resulting equations by $-\nabla\Delta u$ and $-\nabla\Delta v$, respectively, we obtain
\begin{align*}
&\frac{1}{2}\frac{d}{dt}\int_{\mathbb R^N}(|\Delta u|^2+|\Delta v|^2)dx+\varepsilon\int_{\mathbb R^N} |\Delta^2 u|^2dx+ \int_{\mathbb R^N} |\nabla\Delta u|^2dx+\int_{\mathbb R^N} |\nabla\Delta v|^2dx
\\
&=\chi\int_{\mathbb R^N}\nabla(\nabla\cdot(u\nabla v))\nabla\Delta u dx-\xi\int_{\mathbb R^N}\nabla(\nabla\cdot(v\nabla u))\nabla\Delta v dx
\\
&\le\frac12\int_{\mathbb R^N} (|\nabla\Delta u|^2+ |\nabla\Delta v|^2)dx+C\chi^2\int_{\mathbb R^N} (|\nabla^2 v|^2|\nabla u|^2+|\nabla v|^2|\nabla^2u|^2+u^2|\nabla\Delta v|^2)dx
\\
&+C\xi^2\int_{\mathbb R^N} (|\nabla^2 v|^2|\nabla u|^2+|\nabla v|^2|\nabla^2u|^2+v^2|\nabla\Delta u|^2)
\\
&\le\frac12\int_{\mathbb R^N} (|\nabla\Delta u|^2+ |\nabla\Delta v|^2)dx
\\
&+C(\xi^2+\chi^2)\left(\|\nabla  u\|_{L^6}^2\|\nabla^2 v\|_{L^3}^2+\|\nabla v\|_{L^6}^2\|\nabla^2 u\|_{L^3}^2
+\|u\|_{L^\infty}^2\|\nabla\Delta v\|_{L^2}^2+\|v\|_{L^\infty}^2\|\nabla\Delta u\|_{L^2}^2\right)
\\
&\le\frac12\int_{\mathbb R^N} (|\nabla\Delta u|^2+ |\nabla\Delta v|^2)dx+C(\xi^2+\chi^2)\left(\|u\|_{H^2}^2+\|v\|_{H^2}^2\right)
\left(\|\nabla^2u\|_{H^1}^2+\|\nabla^2 v\|_{H^1}^2\right).
\end{align*}
The latter inequality implies that
\begin{align}
&\frac{d}{dt}\int_{\mathbb R^N}(|\Delta u|^2+|\Delta v|^2)dx+2\varepsilon\int_{\mathbb R^N} |\Delta^2 u|^2dx+ \int_{\mathbb R^N} |\nabla\Delta u|^2dx+\int_{\mathbb R^N} |\nabla\Delta v|^2dx
\nonumber
\\
\label{3-8}
&\le C(\xi^2+\chi^2)\left(\|u\|_{H^2}^2+\|v\|_{H^2}^2\right)
\left(\|\nabla^2u\|_{H^1}^2+\|\nabla^2 v\|_{H^1}^2\right).
\end{align}
Now combining the inequalities \eqref{3-4}, \eqref{3-7} and \eqref{3-8}, we then arrive at
the desired inequality \eqref{3-3} readily. \hfill $\Box$

\begin{lemma}\label{lem3-2}
Assume that $\mathcal F(u, \sigma)=u$. Let
$$
\phi(t)=\|u(\cdot, t)\|_{H^2}^{2}+\|v(\cdot, t)\|_{H^2}^{2},
$$
$$
h(t)=\|\nabla u(\cdot,t)\|_{H^2}^2+\|\nabla v(\cdot,t)\|_{H^2}^2.
$$
Then if $\phi(0)\le \frac1{2C^*(\chi^2+\xi^2)}$ with $C^*$ given in Lemma \ref{lem3-1},  we have
$$
\phi(t)\le \phi(0),
$$
$$
2\varepsilon\int_0^\infty \|\Delta u(\cdot,t)\|_{H^2}^2 dt+\frac12 \int_0^\infty h(t)dt\le \phi(0).
$$
\end{lemma}

{\it\bfseries Proof.} From Lemma \ref{lem3-1}, it follows that
$$
\phi'(t)+2\varepsilon\|\Delta u(\cdot,t)\|_{H^2}^2+\Big(1-C^*(\chi^2+\xi^2)\phi(t)\Big)h(t)\le 0.
$$
It is observed that  whenever
$$
\phi(0)\le \frac1{2C^*(\chi^2+\xi^2)},
$$
we have
$$
\phi(t)\le \phi(0)\le \frac1{2C^*(\chi^2+\xi^2)} \quad \text{for all} \ t>0.
$$
and then derive
$$
2\varepsilon \int_0^\infty \|\Delta u(\cdot,t)\|_{H^2}^2 dt+\frac12 \int_0^\infty h(t)dt\le \phi(0),
$$
and thereby complete the proof of this lemma.  \hfill $\Box$

\begin{lemma}\label{lem3-3}
Suppose that $\mathcal F(u, \sigma)=u$ and $\phi(0)\le \frac1{2C^*(\chi^2+\xi^2)}$.
Then we have
$$
\int_0^\infty\int_\Omega (|u_t|^2+|v_t|^2)dxdt\le C,
$$
where $C$ is independent of $\varepsilon$, which depends only on $\Omega$, $u_0$, $v_0$, $\chi$,  and $\xi$.
\end{lemma}

{\it\bfseries Proof.} Multiplying the first equation of \eqref{3-2} by $u_t$, integrating it over $\Omega$ then yields
\begin{align*}
&\frac12\frac{d}{dt}\int_{\mathbb R^N}(\varepsilon|\Delta u|^2+|\nabla u|^2)dx+\int_{\mathbb R^N}|u_t|^2dx
\\
&=-\chi\int_{\mathbb R^N}\nabla\cdot\Big(u\nabla v\Big)u_tdx
\\
&\le \frac12\int_{\mathbb R^N}|u_t|^2dx+C\chi^2\|u\|_{H^2}^2\|\nabla v\|_{H^1}^2.
\end{align*}
Using Lemma \ref{lem3-2} and by a  direct integration, we obtain that
\begin{align*}
\int_0^\infty\int_{\mathbb R^N}|u_t|^2dxdt\le & 2C\chi^2\sup_{t>0}\|u(\cdot, t)\|_{H^2}^2\int_0^\infty\|\nabla v\|_{H^1}^2dt
+\int_{\mathbb R^N}(\varepsilon|\Delta u_0|^2+|\nabla u_0|^2)dx
\\
\le \tilde C.
\end{align*}
Multiplying the second equation of \eqref{3-2} by $v_t$,  integrating the result over $\Omega$ then yields
\begin{align*}
&\frac12\frac{d}{dt}\int_{\mathbb R^N}|\nabla v|^2 dx+\int_{\mathbb R^N}|v_t|^2dx
\\
&=\xi\int_{\mathbb R^N}\nabla\cdot\Big(v\nabla u\Big)v_tdx
\\
&\le \frac12\int_{\mathbb R^N}|v_t|^2dx+C\xi^2\|v\|_{H^2}^2\|\nabla u\|_{H^1}^2.
\end{align*}
From Lemma \ref{lem3-2} we infer that
\begin{align*}
\int_0^\infty\int_{\mathbb R^N}|v_t|^2dxdt
\le 2C\xi^2\sup_{t>0}\|v(\cdot, t)\|_{H^2}^2\int_0^\infty\|\nabla u\|_{H^1}^2dt
+\int_{\mathbb R^N} |\nabla v_0|^2 dx
\le  \hat C,
\end{align*}
and thereby completes the proof.

Furthermore by Lemma \ref{aplem3-3}, we can also conclude that
\begin{lemma}\label{lem3-4}
Assume that $\mathcal F(u, \sigma)=u$, $\phi(0)\le \frac1{2C^*(\chi^2+\xi^2)}$ with  $C^*$  given as in Lemma \ref{lem3-1}. If $|x|^ru_{\varepsilon 0}(x)\in L^2(\mathbb R^N)$ for some a positive constant $r\le 1 $, then  $|x|^ru\in L^\infty_{loc}((0, \infty); L^2(\mathbb R^N))$.
\end{lemma}

\medskip

This lemma allows  us to employ  the compactness result in \cite{ZS}  to demonstrate the existence of strong solutions to system \eqref{1.6}.

\medskip

{\it\bfseries Proof of Theorem \ref{thm-3}.} By the Leray-Schauder fixed point theorem, and Lemmas \ref{lem3-2}--\ref{lem3-4}, the problem
\eqref{3-1} possesses the strong solution $(u_\varepsilon, v_\varepsilon)$ for all $\varepsilon\in [0,1]$, such that $u_\varepsilon\in L^\infty(\mathbb R^+; H^2(\mathbb R^N))\cap L^2(\mathbb R^+; H^4(\mathbb R^N))$, $\partial_tu_\varepsilon\in L^2(\mathbb R^+; L^2(\mathbb R^N))$, $v_\varepsilon\in L^\infty(\mathbb R^+; H^2(\mathbb R^N))\cap L^2(\mathbb R^+; H^3(\mathbb R^N))$, $\partial_tv_\varepsilon\in L^2(\mathbb R^+; L^2(\mathbb R^N))$.
Rechecking the proof of Lemma \ref{lem3-2}, and noticing that
$$
\|\varepsilon\Delta ^2u_{\varepsilon}\|_{L^2(Q)}^2\le C \varepsilon,
$$
$$
\sup_{0<t<T}\int_{\mathbb R^N}(1+|x|^2)^ru^2 dx\le C_T,
$$
letting $\varepsilon \to 0$ (taking its subsequence if necessary),  and using the compact results in \cite{ZS}, for any $T>0$, we obtain
\begin{align*}
&\varepsilon\Delta ^2u_{\varepsilon}\to 0 \ \text{in } \ L^2((0, T); H^2(\mathbb R^N));
\\
&(u_{\varepsilon}, v_{\varepsilon}) \rightharpoonup (u, v) \ \text{in } \ L^2((0, T); H^3(\mathbb R^N));
\\
&(\partial_tu_{\varepsilon}, \partial_tv_{\varepsilon}) \rightharpoonup (\partial_tu, \partial_tv) \ \text{in } \ L^2((0, T); L^2(\mathbb R^N));
\\
&(u_{\varepsilon}, v_{\varepsilon}) \stackrel{*}{\rightharpoonup} (u, v) \ \text{in } \ L^\infty((0, T); H^2(\mathbb R^N));
\\
&(u_{\varepsilon}, v_{\varepsilon}) \rightarrow (u, v) \ \text{in $L^2((0, T); H^2(\mathbb R^N))$};
\\
&(u_{\varepsilon}, v_{\varepsilon}) \rightarrow (u, v) \ \text{uniformly in $Q_T$}.
\end{align*}
Utilizing the respective estimations derived in Lemma \ref{lem3-2}--Lemma \ref{lem3-3}, the limit functions $(u, v)$
with $u, v\in L^\infty(\mathbb R^+; H^2(\mathbb R^N))\cap L^2(\mathbb R^+; H^3(\mathbb R^N))$, $\partial_tu, \partial_tv \in L^2(\mathbb R^+; L^2(\mathbb R^N))$ can be  identified to be the strong solution of \eqref{1.6}. In addition, noticing that $u_0, v_0\ge 0$, then by strong maximum principle, $u(x,t)>0$, $v(x,t)>0$ for $t>0, x\in \mathbb R^N$.

Next, we show the uniqueness of the strong solution of \eqref{1.6}. Indeed, supposed that $(u_1, v_1)$, $(u_2, v_2)$ are the global strong solutions thereof.
Denoting $\hat u=u_1-u_2$, $\hat v=v_1-v_2$, then by a direct calculation, we derive
\begin{align*}
&\frac12\frac{d}{dt}\int_{\mathbb R^N}|\hat u|^2dx+\int_{\mathbb R^N}|\nabla\hat u|^2dx=\chi\int_0^t\int_{\mathbb R^N}(\hat u\nabla v_1\nabla\hat u
+u_2\nabla\hat v\nabla\hat u) dx,
\\
&\frac12\frac{d}{dt}\int_{\mathbb R^N}|\hat v|^2dx+\int_{\mathbb R^N}|\nabla\hat v|^2dx=-\xi\int_0^t\int_{\mathbb R^N}(\hat v\nabla u_1\nabla\hat v
+ v_2\nabla\hat u\nabla\hat v) dx.
\end{align*}
Combining the above two equalities, we further have
\begin{align*}
&\frac12\frac{d}{dt}\int_{\mathbb R^N}\left(|\hat u|^2+|\hat v|^2\right)dx+\int_{\mathbb R^N}\left(|\nabla\hat u|^2+|\nabla\hat v|^2\right)dx
\\
&=\chi\int_{\mathbb R^N}(\hat u\nabla v_1\nabla\hat u
+u_2\nabla\hat v\nabla\hat u) dx-\xi\int_0^t\int_{\mathbb R^N}(\hat v\nabla u_1\nabla\hat v
+ v_2\nabla\hat u\nabla\hat v) dx
\\
&\le \chi \|\hat u\|_{L^3}\|\nabla v_1\|_{L^6}\|\nabla\hat u\|_{L^2}+\xi\|\hat v\|_{L^3}\|\nabla u_1\|_{L^6}\|\nabla\hat v\|_{L^2}
+\chi\|u_2\|_{H^2}\|\nabla\hat v\|_{L^2}\|\nabla\hat u\|_{L^2}
+\xi \|v_2\|_{H^2}\|\nabla\hat v\|_{L^2}\|\nabla\hat u\|_{L^2}
\\
&\le C\chi\|\hat u\|_{L^2}^{\frac{6-N}6}\|\nabla\hat u\|_{L^2}^{\frac N6}\|\nabla v_1\|_{H^1}\|\nabla\hat u\|_{L^2}+\xi\|\hat v\|_{L^2}^{\frac{6-N}6}\|\nabla\hat v\|_{L^2}^{\frac N6}\|\nabla u_1\|_{H^1}\|\nabla\hat v\|_{L^2}
\\
&\quad+\chi\|u_2\|_{H^2}\|\nabla\hat v\|_{L^2}\|\nabla\hat u\|_{L^2}
+\xi \|v_2\|_{H^2}\|\nabla\hat v\|_{L^2}\|\nabla\hat u\|_{L^2}
\\
&\le \frac12 \left(\|\nabla\hat u\|_{L^2}^2+\|\nabla\hat v\|_{L^2}^2\right)+C\left(\|\hat u\|_{L^2}^2+\|\hat v\|_{L^2}^2\right)
+\chi^2\|u_2\|_{H^2}^2\|\nabla\hat v\|_{L^2}^2+
\xi^2\|v_2\|_{H^2}^2\|\nabla\hat u\|_{L^2}^2,
\end{align*}
which implies that
\begin{align*}
&\frac12\frac{d}{dt}\left(\|\hat u\|_{L^2}^2+\|\hat v\|_{L^2}^2\right)-C\left(\|\hat u\|_{L^2}^2+\|\hat v\|_{L^2}^2\right)
\\
&+\left(\frac12-(\chi^2+\xi^2)(\|u_2\|_{H^2}^2+\|v_2\|_{H^2}^2)\right)\int_{\mathbb R^N}\left(|\nabla\hat u|^2+|\nabla\hat v|^2\right)dx\le 0.
\end{align*}
It is observed that $\Phi(t)\le  \frac1{2C^*(\chi^2+\xi^2)}$ with $C^*>1$, and hence
$$
\frac12\frac{d}{dt}\left(\|\hat u\|_{L^2}^2+\|\hat v\|_{L^2}^2\right)-C\left(\|\hat u\|_{L^2}^2+\|\hat v\|_{L^2}^2\right)\le 0.
$$
Therefore,
$$
\|\hat u(\cdot, t)\|_{L^2}^2+\|\hat v(\cdot, t)\|_{L^2}^2\equiv 0,
$$
since $\|\hat u(\cdot, 0)\|_{L^2}^2+\|\hat v(\cdot, 0)\|_{L^2}^2=0$, which thereby completes the proof of the uniqueness of strong solutions.

At last, we pay our attention to the estimate of decay rate claimed in Theorem 1.3.
From \eqref{3-3}, we see that
$$
\frac{d}{dt}\phi(t)+\frac12 h(t)\le 0.
$$
On the other hand, noting that $u, v>0$ and
$$
\int_{\mathbb R^N}u(x,t)dx\equiv \int_{\mathbb R^N}u_0(x)dx=\overline u_0, \quad \int_{\mathbb R^N}v(x,t)dx\equiv \int_{\mathbb R^N}v_0(x)dx=\overline v_0,
$$
and by the Gagliardo-Nirenberg inequality, we infer that
$$
\|u\|_{L^2}^2\le C |\overline u_0|^{\frac{4}{N+2}}\|\nabla u\|_{L^2}^{\frac{2N}{N+2}}, \quad \|v\|_{L^2}^2\le C|\overline v_0|^{\frac{4}{N+2}}\|\nabla u\|_{L^2}^{\frac{2N}{N+2}}.
$$
Therefore, we have
$$
\phi'(t)+A^*\phi(t)^{\frac{N+2}{N}}\le 0,
$$
for some constant $A^*$, and then
$$
\phi(t)\le \left(\Phi^{-\frac 2N}(0)+\frac{2}{N}A^*t\right)^{-\frac{N}2},
$$
which thus completes the proof of Theorem \ref{thm-3} readily. \hfill $\Box$

\begin{appendices}
\section{Appendix}
Consider
\begin{align}
\label{ap-2}\left\{
\begin{aligned}
&v_t-\Delta v=\xi\nabla\cdot(v\nabla\tilde u) && \text{in}\  \mathbb R^N\times \mathbb R^+,
\\
&v(x,0)=\sigma v_{0}(x)\ge 0, &&  x\in\mathbb R^N,
\end{aligned}\right.
\end{align}
where  $\sigma\in[0, 1]$, $v_{0}(x)\in H^2(\mathbb R^N)\cap L^1(\mathbb R^N)$ with $N\le 3$.

\begin{proposition}\label{appro3-1}
Assume that $\tilde u\in L^\infty((0, T); H^2(\mathbb R^N))\cap L^2((0, T); H^3(\mathbb R^N))$. Then \eqref{ap-2} admits a strong solution $v$ with $v\ge 0$, $v\in L^\infty((0, T); H^2(\mathbb R^N)\cap L^1(\mathbb R^N))\cap L^2((0, T); H^3(\mathbb R^N))$, and $v_t \in L^2((0, T); L^2(\mathbb R^N))$.
\end{proposition}

{\it\bfseries Proof.} To prove this lemma, we use bounded domains $B_{R_n}$ to approximate $\mathbb R^N$, where
$B_{R_n}$ denotes the ball centered at $0$ with radius  $R_n>1$, and $R_n\to\infty$ as $n\to\infty$.
In what follows, we denote
$$
Q_n:=B_{R_n}\times (0,\infty); \quad Q_{n, T}:=B_{R_n}\times (0, T);
$$
$\mathring{C}^{\infty}(Q_{n, T})$ denotes the set of functions that are infinitely differentiable and vanish near the lateral boundary of $Q_{n, T}$.

Let $\{v_{n 0}\}_{n=1}^{\infty}$ with $0\le v_{n 0}\in C_0^{\infty}(B_{R_n})$ be the smooth approximation of $v_0$,
$
\{\tilde u_n\}_{n=1}^{\infty}$ with $\tilde u_n\in \mathring{C}^{\infty}(Q_{n, T})$ be the smooth approximation of $\tilde u$,
such that
\begin{align}\label{A2}
\|v_{n,0}\|_{H^2(B_{R_n})}+\|v_{n,0}\|_{L^1(B_{R_n})}\le C,
\end{align}
\begin{align}\label{A3}
\sup_{0<t<T}\|\tilde u_n(\cdot, t)\|_{H^2(B_{R_n})}^2+\int_0^T\|\tilde u_n(\cdot, t)\|_{H^3(B_{R_n})}^2dt\le C,
\end{align}
with $C$ is independent of $n$.
By classical theory of linear parabolic equations, the following problem
\begin{align}\label{ap-1}
\left\{
\begin{aligned}
&v_t-\Delta v=\xi\nabla\cdot(v\nabla\tilde u_n), && x\in B_{R_n}, t>0
\\
&\left.\frac{\partial v}{\partial \nu}\right|_{\partial B_{R_n}}=0
\\
&v(x,0)=\sigma v_{n 0}(x), &&  x \in B_{R_n}
\end{aligned}\right.
\end{align}
admits a unique nonnegative classical solution $v_n\in C^{\infty}(\overline Q_{n,T})$. Therefore, to show the existence of solutions claimed in Proposition \ref{appro3-1}, we concentrate on  derivation of the following a priori estimates on $v_n$
$$
\sup_{0<t<T}\left(\|v_{n}(t)\|_{L^1(B_{R_n})}+\|v_n(t)\|_{H^2(B_{R_n})}^2\right)+ \int_0^T\|v_n(t)\|_{H^3(B_{R_n})}^2dt\le C
$$
with $C$  independent of $n$.

Firstly, it is easy to see that $v_n$  follows the law of conservation of mass, that is
$$
\|v_n(\cdot, t)\|_{L^1(B_{R_n})}\equiv \|v_{n,0}\|_{L^1(B_{R_n})}.
$$
Noticing that
$$
\partial_t v_n-\Delta v_n+v_n=\xi\nabla\cdot(v_n\nabla\tilde u_n)+v_n,
$$
by Neumann heat semigroup theory,  we obtain that
\begin{align*}
&\|v_n(\cdot, t)\|_{L^\infty(B_{R_n})}\le   e^{-t}\|v_{n, 0}\|_{L^\infty(B_{R_n})}
\\
&+\int_0^t e^{-(t-s)}(t-s)^{-\frac{N}{8}-\frac12}\|v_n\nabla\tilde u_n\|_{L^4(B_{R_n})}ds
+\int_0^t e^{-(t-s)}(t-s)^{-\frac{N}{4}}\|v_n\|_{L^2(B_{R_n})}ds
\\
\le & e^{-t}\|v_{n, 0}\|_{L^\infty(B_{R_n})}+\int_0^t e^{-(t-s)}(t-s)^{-\frac{N}{8}-\frac12}\|v_n\|_{L^{12}(B_{R_n})}\|\nabla\tilde u_n\|_{L^6(B_{R_n})}ds
+\int_0^t e^{-(t-s)}(t-s)^{-\frac{N}{4}}\|v_n\|_{L^2(B_{R_n})}ds
\\
\le & e^{-t}\|v_{n, 0}\|_{L^\infty(B_{R_n})}+\sup_{t>0}\{\|v_n(\cdot, t)\|_{L^{12}(B_{R_n})}\|\tilde u_n(\cdot, t)\|_{H^2(B_{R_n})}\}\int_0^t e^{-s}s^{-\frac{N}{8}-\frac12}ds
\\
&+\sup_{t>0}\|v_n(\cdot, t)\|_{L^2(B_{R_n})}\int_0^t e^{-s}s^{-\frac{N}{4}}ds
\\
\le & e^{-t}\|v_{n, 0}\|_{L^\infty(B_{R_n})}+C\sup_{t>0}\|v_n(\cdot, t)\|_{L^{\infty} (B_{R_n})}^{\frac{11}{12}}\|v_n(\cdot, t)\|_{L^{1}(B_{R_n})}^{\frac{1}{12}}
+\sup_{t>0}\|v_n(\cdot, t)\|_{L^{\infty}(B_{R_n})}^{\frac{1}{2}}\|v_n(\cdot, t)\|_{L^{1}(B_{R_n})}^{\frac{1}{2}}.
\end{align*}
Due to $H^2\hookrightarrow L^\infty$ for $N<4$,   $\|v_{n, 0}\|_{L^\infty}$ is bounded by \eqref{A2}. Hence from the above inequality, we infer that
\begin{align}\label{ap-4}
\sup_{t>0}\|v_n(\cdot, t)\|_{L^\infty(B_{R_n})}\le C,
\end{align}
where $C$ is independent of $n$, and thereby
\begin{align}\label{ap-3}
\sup_{t>0}\|v_n(\cdot, t)\|_{L^{r}(B_{R_n})}\le C, \quad \text{for any $1\le r\le \infty$}
\end{align}
since $v_n(\cdot,t)\in L^1(B_{R_n})$, and $C$ is independent of $n$.

Recalling Lemma \ref{GN2-2}, we see that for any $u\in W^{m,r}(B_R)$,
$$
\|D^ju\|_{L^p(B_R)}\le C_1\|D^mu\|_{L^r(B_R)}^\alpha\|u\|^{1-\alpha}_{L^q(B_R)}+C_2R^{\frac Np-\frac Ns-j}\|u\|_{L^s(B_R)}.
$$
From the above inequality, we observe that as long as $\frac Np<\frac Ns+j$, we can replace $R$ with $1$ in this inequality when $R > 1$.
{\bf This ensures that the embedding constants for the Gagliardo-Nirenberg inequality on all bounded domains do not depend on $R_n$ since $R_n>1$.}

Note that
\begin{equation}\label{A.7}
\|\nabla\tilde u_n\|_{L^6(B_{R_n})}^2\le C\|\nabla\tilde u_n\|_{L^2(B_{R_n})}^{\frac{6-2N}3}\|\Delta\tilde u_n\|_{L^2(B_{R_n})}^{\frac{2N}3}+CR_n^{-\frac{2N}3}\|\nabla\tilde u_n\|_{L^2(B_{R_n})}^2
\le C\|\nabla\tilde u_n\|_{H^1(B_{R_n})}^2.
\end{equation}
Multiplying the equation in \eqref{ap-2} by $-\Delta v_n$, and using \eqref{A.7} and Lemma \ref{GN2-2}, we obtain
\begin{align*}
&\frac{1}{2}\frac{d}{dt}\int_{B_{R_n}}|\nabla v_n|^2dx+ \int_{B_{R_n}} |\Delta v_n|^2dx=-\xi\int_{B_{R_n}} \Delta v_n\nabla\cdot(v_n\nabla\tilde u_n) dx
\\
&\le\frac14\int_{B_{R_n}} |\Delta v_n|^2dx+2\xi^2 \int_{B_{R_n}} (v_n^2|\Delta\tilde u_n|^2+|\nabla v_n\nabla\tilde u_n|^2)dx
\\
&\le\frac14\int_{B_{R_n}} |\Delta v_n|^2dx+2\xi^2\|v_n\|_{L^\infty(B_{R_n})}^2\|\Delta\tilde u_n\|_{L^2(B_{R_n})}^2+2\xi^2\|\nabla v_n\|_{L^3(B_{R_n})}^2\|\nabla\tilde u_n\|_{L^6(B_{R_n})}^2
\\
&\le\frac14\int_{B_{R_n}} |\Delta v_n|^2dx+2\xi^2\|v_n\|_{L^\infty(B_{R_n})}^2\|\Delta\tilde u_n\|_{L^2(B_{R_n})}^2
\\
&+C\xi^2\left(\|\Delta v_n\|_{L^2(B_{R_n})}^{\frac{N}3}\|v_n\|_{L^2(B_{R_n})}^{2-\frac{N}3}+\|v_n\|_{L^2(B_{R_n})}^{2}\right)\|\nabla\tilde u_n\|_{H^1(B_{R_n})}^2
\\
&\le \frac12\int_{B_{R_n}} |\Delta v_n|^2dx+\tilde C,
\end{align*}
which implies that
\begin{equation}\label{ap-5}
\sup_{0<t<T}\int_{B_{R_n}}|\nabla v_n|^2dx+ \int_0^T\int_{B_{R_n}} |\Delta v_n|^2dxdt\le C
\end{equation}
with some constant $C$  independent of $n$.

In addition, testing  the equation of \eqref{ap-2} by $\partial_tv_n$, similar to the proof of \eqref{ap-5}, we also have
\begin{equation}\label{ap-6}
\int_0^T\int_{B_{R_n}} |\partial_t v_n|^2dxdt\le C.
\end{equation}
By Lemma \ref{GN2-2}, we obtain the following inequalities,
\begin{align}
\label{ap-8}
&\|\nabla^2 v_n\|_{L^3}^2\le C\|v_n\|_{L^\infty}^{\frac23}\|\nabla\Delta v_n\|_{L^2}^{\frac43}+CR_n^{\frac N3-2}\|v_n\|_{L^\infty}^2
\le \tilde C\left(\|\nabla\Delta v_n\|_{L^2}^{\frac43}+1\right).
\end{align}
Applying $\nabla$ to the equation of \eqref{ap-2}, and multiplying the resultant equation by $-\nabla\Delta v$,  we use \eqref{ap-4}, \eqref{A.7} and \eqref{ap-8}  to get
\begin{align*}
&\frac{1}{2}\frac{d}{dt}\int_{B_{R_n}}|\Delta v_n|^2dx+ \int_{B_{R_n}} |\nabla\Delta v_n|^2dx=-\xi\int_{B_{R_n}}\nabla(\nabla\cdot(v_n\nabla\tilde u_n))\nabla\Delta v_n dx
\\
\le &\frac14\int_{B_{R_n}} |\nabla\Delta v_n|^2dx+C_1\int_{B_{R_n}} (|\nabla^2 v_n|^2|\nabla\tilde u_n|^2+|\nabla v_n|^2|\nabla^2\tilde u_n|^2+v_n^2|\nabla\Delta\tilde u_n|^2)dx
\\
\le&\frac14\int_{B_{R_n}} |\nabla\Delta v_n|^2dx+C_1\Big(\|\nabla\tilde u_n\|_{L^6(B_{R_n})}^2\|\nabla^2 v_n\|_{L^3(B_{R_n})}^2+\|\nabla^2\tilde u_n\|_{L^6(B_{R_n})}^2\|\nabla v_n\|_{L^3(B_{R_n})}^2
\\
&+\|v_n\|_{L^\infty(B_{R_n})}^2\|\nabla\Delta\tilde u_n\|_{L^2(B_{R_n})}^2\Big)
\\
\le & \frac14\int_{B_{R_n}} |\nabla\Delta v_n|^2dx+C_2\Big(\|\nabla\tilde u_n\|_{H^1(B_{R_n})}^2(\|\nabla\Delta v_n\|_{L^2(B_{R_n})}^{\frac{4}3}+1)
\\
&+\|\nabla^2\tilde u_n\|_{H^1(B_{R_n})}^2(\|\nabla v_n\|_{L^2(B_{R_n})}^2+\|\Delta v_n\|_{L^2(B_{R_n})}^2)+
\|\nabla\Delta\tilde u_n\|_{L^2(B_{R_n})}^{2}\Big)
\\
\le &\frac12\int_{B_{R_n}} |\nabla\Delta v_n|^2dx+C_3\left(\|\nabla^2\tilde u_n\|_{H^1(B_{R_n})}^2+\|\nabla^2\tilde u_n\|_{H^1(B_{R_n})}^2\|\Delta v_n\|_{L^2(B_{R_n})}^2
+1\right),
\end{align*}
which implies that
\begin{equation}\label{ap-9}
\sup_{0<t<T}\int_{B_{R_n}}|\Delta v_n|^2dx+ \int_0^T\int_{B_{R_n}} |\nabla\Delta v_n|^2dxdt\le C,
\end{equation}
where $C>0$ is independent of $n$. According to \eqref{ap-4}, \eqref{ap-5}, \eqref{ap-6} and \eqref{ap-9} and by the Cantor diagonal method, one can conclude that there exists a subsequence of $\{v_n\}_{n=1}^\infty$, denoted by itself for  simplicity,  and a nonnegative function  $v\in L^\infty((0, T); H^2(\mathbb R^N)\cap L^1(\mathbb R^N))\cap L^2((0, T); H^3(\mathbb R^N))$, and $v_t \in L^2((0, T); L^2(\mathbb R^N))$ such that
\begin{align*}
&v_n   \stackrel{\star}{\rightharpoonup} v  \quad \text{in}\ L^{\infty}((0, T); H^2_{loc}(\mathbb R^N)\cap L^1_{loc}(\mathbb R^N));
\\
&v_n \rightharpoonup v   \quad \text{in}\
L^2((0,T); H^3_{loc}(\mathbb R^N));\\
&\partial_t v_n \rightharpoonup v_t \quad \ \text{in} \ L^2((0, T); L^2_{loc}(\mathbb R^N));
\\
& v_n\rightarrow v \quad \ \text{in}\  L^2((0, T); H^2_{loc}(\mathbb R^N)),
\end{align*}
and $v$ is exactly a strong solution of \eqref{ap-2}.
 \hfill $\Box$

\medskip

Now for  the function $v(x,t)$  obtained  in Proposition \ref{appro3-1}, we consider
\begin{align}
\label{ap-10}\left\{
\begin{aligned}
&u_t+\varepsilon\Delta^2u-\Delta u=-\chi\nabla\cdot(u\nabla v) && \text{in}\  \mathbb R^N\times \mathbb R^+,
\\
&u(x,0)=\sigma u_{0}(x), &&  x\in\mathbb R^N,
\end{aligned}\right.
\end{align}
with $xu_{0}(x)\in L^2(\mathbb R^N)$, $u_{0}(x)\in H^3(\mathbb R^N)$.

\begin{proposition}\label{appro3-2}
Assume that $v\in L^\infty((0, T); H^2(\mathbb R^N))\cap L^2((0, T); H^3(\mathbb R^N))$. Then \eqref{ap-10}
admits a strong solution   $u\in L^\infty((0, T); H^3(\mathbb R^N))\cap L^2((0, T); H^5(\mathbb R^N))$ and
$u_t\in L^2(Q_T)$.
\end{proposition}

From Lemma \ref{GN2-2} and the proof of Proposition \ref{appro3-2}, we can observe that the embedding constants in either the Gagliardo-Nirenberg inequality or the Sobolev inequality (see for example \eqref{A.7}, \eqref{ap-8} etc.) can be selected in a manner that is independent of
$R_n$. Consequently, the subsequent verification of this proposition adheres to this principle. For simplicity, in the subsequent proofs, we will not provide detailed explanations for each inequality as we did in the previous proposition.

\medskip

{\it\bfseries Proof.} As the proof in Proposition \ref{appro3-1}, in order to obtain the a prior estimates for
solutions  of approximate system of \eqref{ap-10}, we consider the  initial boundary problem of the form
\begin{align}\label{app-9}
\left\{
\begin{aligned}
&u_t+\varepsilon\Delta^2u-\Delta u=-\chi\nabla\cdot(u\nabla\tilde v_n), && x\in B_{R_n}, t>0,
\\
&\left.\frac{\partial u}{\partial \nu}\right|_{\partial B_{R_n}}=0,  \left.\frac{\partial \Delta u}{\partial \nu}\right|_{\partial B_{R_n}}=0,
\\
&u(x,0)=\sigma u_{n 0}(x), \quad x \in B_{R_n}, &&
\end{aligned}\right.
\end{align}
where $\{u_{n 0}\}_{n=1}^{\infty}$ with $u_{n 0}\in C_0^{\infty}(B_{R_n})$ are the smooth approximation of $u_0$,
$\{\tilde v_n\}_{n=1}^{\infty}$ with $\tilde v_n\in \mathring{C}^{\infty}(Q_{n, T})$ are the smooth approximation of $v$,
such that
$$
\|u_{n,0}\|_{H^3(B_{R_n})}\le C,
$$
$$
\sup_{0<t<T}\|\tilde v_n(\cdot, t)\|_{H^2(B_{R_n})}^2+\int_0^T\|\tilde v_n(\cdot, t)\|_{H^3(B_{R_n})}^2dt\le C,
$$
with $C$ independent of $n$.
By classical theory of linear parabolic equations, the problem \eqref{app-9} admits a unique  classical solution $u_n\in C^{\infty}(\overline Q_{n,T})$. Therefore, to show Proposition \ref{appro3-2}, we only need to derive the following a priori estimate
$$
\sup_{0<t<T}\left(\|u_n(t)\|_{H^3(B_{R_n})}^2\right)+ \int_0^T\left(\|u_n(t)\|_{H^5(B_{R_n})}^2+\|\partial_t u_n(t)\|_{L^2(B_{R_n})}^2\right)dt\le C,
$$
with $C$  independent of $n$.

Multiplying  the equation of  \eqref{app-9} by $u_n$, then integrating the resultant equation over $B_{R_n}$, and using Sobolev embedding
inequality, we obtain
\begin{align*}
&\frac{1}{2}\frac{d}{dt}\int_{B_{R_n}}|u_n|^{2}dx+\varepsilon\int_{B_{R_n}} |\Delta u_n|^2dx+ \int_{B_{R_n}} |\nabla u_n|^2dx
\\
&= -\frac{\chi}{2}\int_{B_{R_n}} u_n^2\Delta \tilde v_n
\\
&\le \frac{\chi}{2}\|\Delta \tilde v_n\|_{L^2(B_{R_n})}\|u_n\|_{L^4(B_{R_n})}^2
\\
&\le \frac12\int_{B_{R_n}} |\nabla u_n|^2dx+C\int_{B_{R_n}}|u_n|^{2}dx,
\end{align*}
which implies that
\begin{align}\label{ap-11}
\sup_{0<t<T}\|u_n(\cdot, t)\|_{L^{2}(B_{R_n})}^2+\int_0^T\int_{B_{R_n}}(\varepsilon|\Delta u_n|^2+|\nabla u_n|^2)dxdt\le C,
\end{align}
where $C$  is independent of $n$.

Testing \eqref{app-9} by $-\Delta u_n$, and integrating the resultant equation over $B_{R_n}$, then the Sobolev embedding
inequality gives us
\begin{align*}
&\frac{1}{2}\frac{d}{dt}\int_{B_{R_n}}|\nabla u_n|^{2}dx+\varepsilon\int_{B_{R_n}} |\nabla\Delta u_n|^2dx+ \int_{B_{R_n}}|\Delta u_n|^2dx
\\
&= -\chi \int_{B_{R_n}} u_n\nabla\tilde v_n\nabla\Delta u_n dx
\\
&\le \frac{\varepsilon}{2}\int_{B_{R_n}}|\nabla\Delta u_n|^2dx+\frac{\chi^2}{2\varepsilon}\|u_n\|_{L^3}^2\|\nabla\tilde v_n\|_{L^6}^2
\\
&\le \frac{\varepsilon}{2}\int_{B_{R_n}}|\nabla\Delta u_n|^2dx+\frac{\chi^2}{2\varepsilon}(C\|u_n\|_{L^2(B_{R_n})}^2+C\|\nabla u_n\|_{L^2(B_{R_n})}^2)\|\nabla\tilde v_n\|_{H^1(B_{R_n})}^2
\\
&\le \frac{\varepsilon}{2}\int_{B_{R_n}}|\nabla\Delta u_n|^2dx+C_\varepsilon(\|\nabla u_n\|_{L^2(B_{R_n})}^2+1),
\end{align*}
which implies that
\begin{align}\label{ap-12}
\sup_{0<t<T}\|\nabla u_n(\cdot, t)\|_{L^{2}(B_{R_n})}^2+\int_0^T\int_{B_{R_n}}(\varepsilon|\nabla\Delta u_n|^2+|\Delta u_n|^2)dxdt\le C,
\end{align}
where $C$ is independent of $n$.
Noticing that
\begin{align*}
\|\nabla \cdot(u_n\nabla\tilde v_n)\|_{L^2(B_{R_n})}^2\le & \|u_n\|_{L^3}^2\|\Delta\tilde v_n\|_{L^6}^2+\|\nabla u_n\|_{L^3}^2\|\nabla\tilde v_n \|_{L^6}^2
\\
\le &C_1(\|u_n\|_{H^1(B_{R_n})}^2\|\Delta\tilde v_n\|_{H^1(B_{R_n})}^2+\|\nabla u_n\|_{L^2(B_{R_n})}^2\|\nabla\tilde v_n \|_{H^1(B_{R_n})}^2+\|\Delta u_n\|_{L^2(B_{R_n})}^2\|\nabla\tilde v_n \|_{H^1(B_{R_n})}^2)
\\
\le & C_2(\|\Delta\tilde v_n\|_{H^1(B_{R_n})}^2+\|\Delta u_n\|_{L^2(B_{R_n})}^2+1),
\end{align*}
recalling that $\tilde v_n\in L^2(B_{R_n})((0, T); H^3(B_{R_n} ))$, and \eqref{ap-12},  we then can see that
\begin{align}\label{ap-13}
\int_0^T\int_{B_{R_n}} |\nabla \cdot(u_n\nabla\tilde v_n)|^2dxdt\le C.
\end{align}
Multiplying the equation in \eqref{app-9} by $\Delta^2 u_n$, and using \eqref{ap-11}, \eqref{ap-12}, we obtain that
\begin{align*}
&\frac{1}{2}\frac{d}{dt}\int_{B_{R_n}}|\Delta u_n|^{2}dx+\varepsilon\int_{B_{R_n}} |\Delta^2 u_n|^2dx+ \int_{B_{R_n}}|\nabla\Delta u_n|^2dx
\\
=& -\chi \int_{B_{R_n}} \nabla \cdot(u_n\nabla\tilde v_n)\Delta^2 u_n dx
\\
\le & \frac{\varepsilon}{2}\int_{B_{R_n}}|\Delta^2 u_n|^2dx+\frac{\chi^2}{\varepsilon}\int_{B_{R_n}} |\nabla \cdot(u_n\nabla\tilde v_n)|^2dxdt.
\end{align*}
From \eqref{ap-13}, we infer that
\begin{align}\label{ap-14}
\sup_{0<t<T}\|\Delta u_n(\cdot, t)\|_{L^{2}(B_{R_n})}^2+\int_0^T\int_{B_{R_n}}(\varepsilon|\Delta^2 u_n|^2+|\nabla\Delta u_n|^2)dxdt\le C,
\end{align}
where $C$ is independent of $n$.

Similarly, by multiplying \eqref{app-9} by $\partial_tu_n$ and using \eqref{ap-13}, we also have
\begin{align}\label{ap-15}
\int_0^T\int_{B_{R_n}} |\partial_tu_n|^2 dxdt\le C,
\end{align}
where $C$ is independent of $n$.

Recalling that $\tilde v_n\in L^\infty((0, T); H^2(B_{R_n} ))\cap L^2((0, T); H^3(B_{R_n} ))$, using \eqref{ap-11}, \eqref{ap-12} and \eqref{ap-14}, we obtain that
\begin{align*}
\|\nabla\nabla \cdot(u_n\nabla\tilde v_n)\|_{L^2(B_{R_n})}^2\le & C_1(\|\nabla^2u_n\|_{L^3(B_{R_n})}^2\|\nabla\tilde v_n\|_{L^6}^2+\|\nabla u_n\|_{L^3(B_{R_n})}^2\|\nabla^2\tilde v_n \|_{L^6(B_{R_n})}^2
+\|u_n\|_{L^\infty(B_{R_n})}^2\|\nabla\Delta\tilde v_n\|_{L^2(B_{R_n})}^2)
\\
\le &C_2(\|u_n\|_{H^3(B_{R_n})}^2\|\tilde v_n\|_{H^2(B_{R_n})}^2+\|u_n\|_{H^2(B_{R_n})}^2\|\tilde v_n\|_{H^3(B_{R_n})}^2)
\\
\le & C_3(\|u_n\|_{H^3(B_{R_n})}^2+\|\tilde v_n\|_{H^3(B_{R_n})}^2+1),
\end{align*}
which implies that
\begin{align}\label{ap-16}
\int_0^T\|\nabla\nabla \cdot(u_n\nabla\tilde v_n)\|_{L^2(B_{R_n})}^2 dt\le C.
\end{align}
Applying $\nabla$ to the  equation of \eqref{app-9}, and multiplying the resultant equation by $\nabla\Delta^2 u$,
we arrive at
\begin{align*}
&\frac{1}{2}\frac{d}{dt}\int_{B_{R_n}}|\nabla\Delta u_n|^{2}dx+\varepsilon\int_{B_{R_n}} |\nabla\Delta^2 u_n|^2dx+ \int_{B_{R_n}}|\Delta^2 u_n|^2dx
\\
&= -\chi \int_{B_{R_n}} \nabla\nabla \cdot(u_n\nabla\tilde v_n)\nabla\Delta^2 u_n dx
\\
&\le \frac{\varepsilon}2\int_{B_{R_n}} |\nabla\Delta^2 u_n|^2dx+\frac{\chi^2}{2\varepsilon}\int_{B_{R_n}} |\nabla\nabla \cdot(u_n\nabla\tilde v_n)|^2dx,
\end{align*}
which implies that
\begin{align}\label{ap-17}
\sup_{0<t<T}\|\nabla\Delta u_n(\cdot, t)\|_{L^{2}(B_{R_n})}^2+\int_0^T\int_{B_{R_n}}(\varepsilon|\nabla\Delta^2 u_n|^2+|\Delta^2 u_n|^2)dxdt\le C,
\end{align}
where $C$ is independent of $n$.
Similar to the proof of Proposition \ref{appro3-1}, we thus complete the proof upon letting $n\to\infty$ in \eqref{ap-17}. \hfill $\Box$

\begin{lemma}\label{aplem3-3}
Assume that $v\in L^\infty((0, T); H^2(\mathbb R^N))\cap L^2((0, T); H^3(\mathbb R^N))$. Let $u\in L^\infty((0, T); H^3(\mathbb R^N))\cap L^2((0, T); H^5(\mathbb R^N))$ be the corresponding strong solution of \eqref{ap-10}
such that
$$
\|u(\cdot, t)\|_{H^2}\le C,
$$
where $C$ is independent of $\varepsilon$. Furthermore, if $(1+|x|^2)^{r/2}u_0 \in L^2(\mathbb R^N)$ for some positive constant $r\le 1 $, then
\begin{align*}
\sup_{0<t<T}\int_{\mathbb R^N}(1+|x|^2)^ru^2 dx+\varepsilon \int_0^T\int_{\mathbb R^N}(1+|x|^2)^r|\Delta u|^2 dxdt+\int_0^T\int_{\mathbb R^N}(1+|x|^2)^r|\nabla u|^2 dxdt\le C_T,
\end{align*}
where $C_T$ is constant depending on $T$ but independent of $\varepsilon$.
\end{lemma}

{\it\bfseries Proof.} For any $R>1$, let $\eta_R(x)\in C_0^{2}(B_{2R})$ be a cut off function on $B_{2R}$ such that $0\le \eta_R(x)\le 1$, $\eta_R(x)=1$ for $x\in B_{R}$, $|\nabla^k\eta_R(x)|\le \frac{C}{R^k}$, where $C$ is a constant independent of $R$.
Let
$$
\varphi(x)=|\eta_R(x)|^q(1+|x|^2)^{r},
$$
where $q$ is a positive constant to be determined below, $r\in (0,1)$.
By a direct calculation, we obtain
\begin{equation}\label{ap-18}
\frac{|\nabla \varphi|^2}{\varphi}\le 8r^2(1+|x|^2)^{r-2}|x|^2\eta_R^q+2q^2(1+|x|^2)^{r}\eta_R^{q-2}|\nabla\eta_R|^2<8+C\eta_R^{q-2}\frac{(1+4R^2)^{r}}{R^2},
\end{equation}
\begin{align}
\frac{|\Delta \varphi|}{\varphi}\le & C\left(1+\eta_R^{q-2}|\nabla\eta_R|^2+\eta_R^{q-4}|\nabla\eta_R|^4(1+|x|^2)^{r}+\eta_R^{q-2}|\Delta\eta_R|^2(1+|x|^2)^{r}\right)\nonumber
\\
\label{ap-19}
\le & C\left(1+\frac 1{R^2}\eta_R^{q-2}+\eta_R^{q-4}\frac{(1+4R^2)^{r}}{R^4}\right).
\end{align}
\begin{equation}\label{ap-20}
|\nabla \varphi|\le 2r\eta_R^q(1+|x|^2)^{r-1}|x|+q\eta_R^{q-1}(1+|x|^2)^{r}|\nabla\eta_R|\le 2\varphi+q\eta_R^{q-1}(1+|x|^2)^{r}|\nabla\eta|.
\end{equation}
In what follows, we take $q=4$. Then
\begin{equation}\label{ap-21}
\frac{|\nabla \varphi|^2}{\varphi}\le C,\quad \frac{|\Delta \varphi|}{\varphi}\le C.
\end{equation}
Multiplying the first equation of \eqref{ap-10} by $\varphi(x)u$, and integrating it over $\mathbb R^N$, and using \eqref{ap-20}, \eqref{ap-21}, we arrive at
\begin{align*}
&\frac12\frac{d}{dt}\int_{B_{2R}}\varphi(x)u^2 dx+\varepsilon\int_{B_{2R}}\varphi(x)|\Delta u|^2 dx+\int_{B_{2R}}\varphi(x)|\nabla u|^2 dx
\\
&=-2\varepsilon\int_{B_{2R}}\nabla\varphi\nabla u\Delta udx-\varepsilon\int_{B_{2R}}u\Delta\varphi\Delta udx-
\int_{B_{2R}}u\nabla\varphi\nabla u dx+\chi\int_{B_{2R}}u\nabla v(\varphi\nabla u+u\nabla\varphi)dx
\\
&\le \frac{\varepsilon}2\int_{B_{2R}}\varphi|\Delta u|^2 dx+4\varepsilon\int_{B_{2R}}\frac{|\nabla\varphi|^2}{\varphi}|\nabla u|^2 dx
+\varepsilon\int_{B_{2R}}\frac{|\Delta\varphi|^2}{\varphi}|u|^2 dx+\frac12\int_{B_{2R}}\varphi|\nabla u|^2 dx+
\int_{B_{2R}}\frac{|\nabla\varphi|^2}{\varphi}|u|^2 dx
\\
&
+\chi^2\|\nabla v\|_{L^\infty}^2\int_{B_{2R}}\varphi u^2dx+2\chi\|\nabla v\|_{L^\infty}\int_{B_{2R}}\varphi u^2 dx
+4\chi\|\nabla v\|_{L^\infty}\int_{B_{2R}}\eta_R^{3}(1+|x|^2)^{r}|\nabla\eta|u^2 dx
\\
&\le \frac{\varepsilon}2\int_{B_{2R}}\varphi(x)|\Delta u|^2 dx+C\varepsilon\int_{B_{2R}}|\nabla u|^2 dx
+C\int_{B_{2R}}|u|^2 dx+\frac12\int_{B_{2R}}\varphi(x)|\nabla u|^2 dx
\\
&
+(\chi^2\|\nabla v\|_{H^2}^2+2\chi\|\nabla v\|_{L^\infty})\int_{B_{2R}}\varphi u^2dx+\chi^2\|\nabla v\|_{H^2}^2\int_{B_{2R}}\varphi u^2dx+
4\int_{B_{2R}} \eta_R^{2}(1+|x|^2)^{r}|\nabla\eta|^2 u^2 dx
\\
&\le \frac{\varepsilon}2\int_{B_{2R}}\varphi(x)|\Delta u|^2 dx+\frac12\int_{B_{2R}}\varphi(x)|\nabla u|^2 dx
+(3\chi^2\|\nabla v\|_{H^2}^2+1)\int_{B_{2R}}\varphi u^2dx
\\
&+C\varepsilon\int_{B_{2R}}|\nabla u|^2 dx
+C\int_{B_{2R}}|u|^2 dx
+4\frac{(1+|4R|^2)^{r}}{R^2}\int_{B_{2R}} u^2 dx
\\
&\le \frac{\varepsilon}2\int_{B_{2R}}\varphi(x)|\Delta u|^2 dx+\frac12\int_{B_{2R}}\varphi(x)|\nabla u|^2 dx
+(3\chi^2\|\nabla v\|_{H^2}^2+1)\int_{B_{2R}}\varphi u^2dx+C.
\end{align*}
That is
\begin{align}
&\frac{d}{dt}\int_{B_{2R}}\varphi(x)u^2 dx+\varepsilon\int_{B_{2R}}\varphi(x)|\Delta u|^2 dx+\int_{B_{2R}}\varphi(x)|\nabla u|^2 dx\nonumber
\\
\label{ap-22}
&\le (6\chi^2\|\nabla v\|_{H^2}^2+2)\int_{B_{2R}}\varphi u^2dx+C,
\end{align}
where $C$ is independent of $R$ and $\varepsilon$. Ignoring the last two terms on the left side of the above inequality, a direct calculation yields
$$
\frac{d}{dt}\left\{\int_{B_{2R}}\varphi u^2 dx \exp\left\{-\int_0^t(6\chi^2\|\nabla v(s)\|_{H^2}^2+2)ds\right\}\right\}
\le C\exp\left\{-\int_0^t(6\chi^2\|\nabla v(s)\|_{H^2}^2+2)ds\right\}.
$$
The integration of the above inequality from $0$ to $t$ gives
\begin{align}
\label{ap-23}
\sup_{0<t<T}\int_{B_{2R}}\varphi u^2 dx\le C,
\end{align}
where $C$ is independent of $R$ and $\varepsilon$.
Recalling \eqref{ap-22}, integrating it from $0$ to $T$, and using \eqref{ap-23} gives
 \begin{align}
&\int_{B_{2R}}\varphi(x)u^2 dx+\varepsilon\int_0^T\int_{B_{2R}}\varphi(x)|\Delta u|^2 dxdt+\int_0^T\int_{B_{2R}}\varphi(x)|\nabla u|^2 dxdt\nonumber
\\
&\le \sup_{0<t<T}\int_{B_{2R}}\varphi u^2dx\int_0^T(6\chi^2\|\nabla v\|_{H^2}^2+2)dt+\int_{B_{2R}}\varphi(x)u_0^2 dx+CT\nonumber
\\
\label{ap-24}
&\le \tilde C,
\end{align}
where $\tilde C$ is independent of $R$ and $\varepsilon$. Letting $R\to\infty$, we achieve  the desired inequality. \hfill $\Box$

\end{appendices}
\vskip3mm

\bf{Acknowledgements}\rm:
Chunhua Jin is supported by the National Natural Science Foundation of China (No.12271186, 12171166). Yifu Wang is supported by the National Natural Science Foundation of China (No.12071030).
\vskip2mm
\bf{Data availability}\rm: Data sharing is not applicable to this article, as no datasets were
created or analyzed as part of the current study.
\vskip2mm
\bf{Conflict of interest}\rm: There is no conflict of interest to this work.

\end{document}